\documentclass[11pt]{article}

\author{Viet Son Pham\thanks{Center for Mathematical Sciences, Technical University of Munich, Boltzmannstra\ss e 3, 85748 Garching, Germany, e-mail: vietson.pham@tum.de}}
\title{L{\'e}vy-driven causal CARMA random fields}
\date{\today}

\usepackage[T1]{fontenc}
\usepackage[utf8]{inputenc}
\usepackage{amsmath}
\usepackage{amsfonts}
\usepackage{amssymb}
\usepackage{amsthm}
\usepackage{amsopn}		
\usepackage{mathrsfs}	
\usepackage{dsfont}
\usepackage{color}
\usepackage{comment}
\usepackage[numbers]{natbib}
\usepackage{longtable}
\usepackage{enumitem}
\setenumerate{leftmargin=*}
\usepackage[linktocpage]{hyperref}
\usepackage{graphicx}
\usepackage{float}

\numberwithin{equation}{section}	
\allowdisplaybreaks[3]	
\newcommand{\bbn}{\mathbb{N}}
\newcommand{\bbz}{\mathbb{Z}}
\newcommand{\bbr}{\mathbb{R}}
\newcommand{\bbc}{\mathbb{C}}
\newcommand{\bbe}{\mathbb{E}}

\newcommand{\bbp}{\mathbb{P}}

\newcommand{\bone}{\mathds 1}

\newcommand{\calb}{{\cal B}}

\newcommand{\cald}{{\cal D}}

\newcommand{\calf}{{\cal F}}

\newcommand{\calo}{{\cal O}}

\newcommand{\cals}{{\cal S}}

\newcommand{\al}{{\alpha}}
\newcommand{\ga}{{\gamma}}

\newcommand{\la}{{\lambda}}
\newcommand{\La}{{\Lambda}}
\newcommand{\si}{{\sigma}}

\newcommand{\om}{{\omega}}

\newcommand{\bi}{\boldsymbol{i}}

\newcommand{\bs}{\textbf{s}}

\newcommand{\frakD}{\mathfrak{D}}
\newcommand{\frakp}{\mathfrak{p}}
\newcommand{\frakq}{\mathfrak{q}}

\newcommand{\bbb}{\mathrm{b}}
\newcommand{\dd}{\mathrm{d}}
\newcommand{\ee}{\mathrm{e}}

\newcommand{\var}{{\mathrm{Var}}}	
\newcommand{\cov}{{\mathrm{Cov}}}	

\newcommand{\halmos}{\quad\hfill $\Box$}	

\newcommand{\Leb}{\mathrm{Leb}}
\newcommand{\loc}{\mathrm{loc}}

\newtheoremstyle{neu}
{11pt}      
{11pt}      
{}                  
{}          
{\bfseries} 
{}          
{1em}  
{\textbf{\thmname{#1}\thmnumber{ #2}\thmnote{ (#3)}}}          

\newtheoremstyle{proof}
{11pt}      
{11pt}      
{}                  
{}          
{\bfseries} 
{}            
{1em}          
{\textbf{\thmname{#1}.}}          

\newtheorem{Theorem}{Theorem}[section]
\newtheorem{Corollary}[Theorem]{Corollary}
\newtheorem{Lemma}[Theorem]{Lemma}
\newtheorem{Proposition}[Theorem]{Proposition}
\theoremstyle{neu}
\newtheorem{Definition}[Theorem]{Definition}
\newtheorem{Example}[Theorem]{Example}
\newtheorem{Remark}[Theorem]{Remark}

\theoremstyle{proof}
\newtheorem{Proof}{Proof}

\begin{document}

\maketitle

\begin{abstract}
	We introduce Lévy-driven causal CARMA random fields on $\bbr^d$, extending the class of CARMA processes. The definition is based on a system of stochastic partial differential equations which generalize the classical state-space representation of CARMA processes. The resulting CARMA model differs fundamentally from the isotropic CARMA random field of Brockwell and Matsuda. We show existence of the model under mild assumptions and examine some of its features including the second-order structure and path properties. In particular, we investigate the sampling behavior and formulate conditions for the causal CARMA random field to be an ARMA random field when sampled on an equidistant lattice. 
\end{abstract}

\vfill

\noindent
\begin{tabbing}
	\> secondary: \,\,\,60G51, 60J75
\end{tabbing}

\vspace{1cm}


\vspace{0.5cm}

\newpage

\section{Introduction}

Originally appearing in \citet{Doob44}, continuous-time autoregressive moving average processes, or CARMA processes in short, are the continuous-time analogs of the well-known ARMA processes (cf. \citet{Brockwell91} for details on ARMA processes). Nowadays, CARMA processes are well-studied objects due to the extensive research in recent years, which partially stems from the successful usage of these processes as stochastic models for irregularly spaced or high-frequency data (cf. the review article \citet{Brockwell14} and the references therein). Applications can be found in turbulence modeling \cite{Brockwell13}, stochastic volatility modeling \cite{BN01a,Brockwell11,Todorov06} and the electricity market \cite{Benth14,Garcia10}, just to name a few.

Given two non-negative integers $q<p$ and real coefficients $a_1,...,a_p,b_0,...,b_{p-1}$ such that $b_q\neq 0$ and $b_i=0$ for $i>q$, the CARMA$(p,q)$ process $(Y(t))_{t\in\bbr}$ is defined as the solution to the observation and state equations
\begin{equation}\label{statespace1}
\begin{aligned}
Y(t)&=b^\top X(t),\quad t\in\bbr,\\
\dd X(t)&=AX(t)\,\dd t+c\,\dd L(t),\quad t\in\bbr,
\end{aligned}
\end{equation}
where $b=(b_0,...,b_{p-1})^\top\in\bbr^p$, $c=(0,...,0,1)^\top\in\bbr^p$ and the matrix $A$ is given by 
\begin{equation*}
A=\begin{pmatrix}
0	&1	&0	&\cdots	&0\\
0	&0	&1	&\cdots	&0\\
\vdots	&\vdots	&\vdots	&\ddots	&\vdots\\
0	&0	&0	&\cdots	&1\\
-a_p	&-a_{p-1}	&-a_{p-2}	&\cdots	&-a_1\\
\end{pmatrix}\in\bbr^{p\times p}
\end{equation*}
if $p>1$, and $A=(-a_1)$ if $p=1$. Here it is assumed that $L$ is a one-dimensional Lévy process, that is, a process with independent and stationary increments, càdlàg sample paths and $L(0)=0$ almost surely (cf. \citet{Sato99} for details on Lévy processes). Equations \eqref{statespace1} can be interpreted as the $p$th-order stochastic differential equation
\begin{equation}
\label{CARMAeq}a(\partial_t)Y(t)=b(\partial_t)\partial_tL(t),\quad t\in\bbr,
\end{equation}
where the polynomials $a(\cdot)$ and $b(\cdot)$ are defined as
\begin{equation}
\label{polynomials} a(z)=z^p+a_1z^{p-1}+\cdots+a_p,\quad\text{and}\quad
b(z)=b_0+b_1z+\cdots+b_{p-1}z^{p-1}.
\end{equation}
In fact, Equation \eqref{CARMAeq} constitutes the continuous-time analog of the well-known ARMA equations, which define the ARMA process in discrete time. However, since the sample paths of a Lévy process are in general not differentiable, the definition of CARMA processes is based on the state-space representation \eqref{statespace1}. Note that $A$ is the companion matrix of the polynomial $a(\cdot)$ and therefore the eigenvalues of $A$ are equal to the roots of $a(\cdot)$. Under the assumptions that $a(\cdot)$ and $b(\cdot)$ have no common roots, the roots of $a(\cdot)$ have strictly negative real parts and the Lévy process $L$ has a finite logarithmic moment, it was shown in \citet[Theorem~3.3]{Brockwell09} that the CARMA equations have a unique strictly stationary solution $Y$ on $\bbr$ with representation
\begin{equation}\label{CARMAproc}
Y(t)=\int_{-\infty}^t b^\top\ee^{A(t-s)}c \,\dd L(s),\quad t\in\bbr.
\end{equation}
Moreover, the CARMA process is a causal function of the driving Lévy process under the assumptions above, i.e., the value of $Y(t)$ depends only on the values of $(L(s))_{s\leq t}$ and is independent of $(L(s))_{s>t}$.

The aim of this article is to extend CARMA processes to multiple parameters in order to obtain a tractable class of random fields indexed by $\bbr^d$, which can be used to model spatial or even temporo-spatial phenomena. A spatial extension has in fact already been introduced in \citet{Brockwell17}. Their isotropic CARMA random field is defined as 
\begin{equation}\label{isoCARMArf}
Y(t)=\int_{\bbr^d}g(t-s)\dd L(s),\quad t\in\bbr^d,
\end{equation}
where the radially symmetric kernel $g$ is given by $g(t)=\sum_{i=1}^p \ee^{\la_i\|t\|}\theta(\la_i)/\phi'(\la_i)$, $t\in\bbr^d$, and the polynomials $\phi(\cdot)$ and $\theta(\cdot)$ have the forms
$\phi(z)=\prod_{i=1}^p (z^2-\la_i^2)$ and $\theta(z)=\prod_{i=1}^q (z^2-\xi_i^2)$ with $\la_i,\xi_i\in\bbc$. Furthermore, each $\la_i$ has a strictly negative real part and $L$ is a Lévy sheet on $\bbr^d$, which is the multi-parameter analog of a Lévy process.
This procedure generates a versatile family of isotropic covariance functions in space, which are neither necessarily monotonically decreasing nor non-negative. However, the CARMA process is classically defined through the state-space representation \eqref{statespace1}, and the connection to these defining equations is unclear in \cite{Brockwell17}. Therefore, we propose a different class of CARMA random fields based on a system of stochastic partial differential equations (SPDEs) that constitutes a generalization of \eqref{statespace1}. We will show that this system has a \emph{mild solution} given by
\begin{equation}\label{CARMArf}
Y(t)=\int_{-\infty}^{t_1}\cdots\int_{-\infty}^{t_d} b^\top\ee^{A_1(t_1-s_1)}\cdots\ee^{A_d(t_d-s_d)}c \,\dd L(s),\quad t=(t_1,...,t_d)\in\bbr^d,
\end{equation}
and define the random field $Y$ in \eqref{CARMArf} as the \emph{causal CARMA random field} on $\bbr^d$, where $A_1,...,A_d$ are companion matrices.

It turns out that many of the commonly known features of CARMA processes can be recovered for this model, including for instance exponentially decaying autocovariance functions and rational spectral densities. Moreover, the autocovariance is in general anisotropic and non-separable. 
The path properties are also similar to those we have in the one-dimensional case. More precisely, there exists a Hölder continuous version under Gaussian noise, and in the presence of jumps, we may use maximal inequalities for multi-parameter martingales in order to show the existence of càdlàg sample paths (see Definition~\ref{cadlag}). Furthermore, \eqref{CARMArf} reduces to \eqref{CARMAproc} if $d=1$ and sampling on an equidistant lattice leads to an ARMA random field under mild conditions. However, the moving average part has in general infinitely many terms in contrast to the one-dimensional case. This is due to the fact that a $(q_1,q_2)$-dependent random field is not always a MA$(q_1,q_2)$ random field (see Definition~\ref{DefARMA}). We examine this issue in Examples~\ref{(1,1)-dependent} and \ref{(1,1)-dependent 2}.

The CARMA random field in this article is causal with respect to the spatial partial order $\leq$ on $\bbr^d$, which is taken componentwise. This quarter-plane-type causality can be interpreted as a directional influence and has been incorporated in several articles in the literature. For example, both \citet{Tjostheim78} and \citet{Drapatz16} consider quarter-plane ARMA models and they refer to applications in econometrics, veterinary epidemiology, geography, geology and image analysis. Furthermore, \cite{Tjostheim78} points out that causal representations exist for a wide class of random fields and \cite{Drapatz16} mentions that statistical inference for such representations is often easier to conduct.

The causal CARMA random field is related to some other classes of random fields. It belongs to the class of \emph{ambit fields} (see e.g. \citet{BN16}), which has applications in biology, finance or turbulence. In particular,  causal CARMA random fields possess ambit sets which are translation invariant and have the form of a \emph{quadrant} if $d=2$, an \emph{octant} if $d=3$, an \emph{orthant} if $d>3$, respectively.
Additionally, we will see that they constitute a parametric submodel of the Volterra-type Ornstein-Uhlenbeck (VOU) processes studied in \citet{Pham18} and they generalize the multi-parameter Ornstein-Uhlenbeck process in \citet{Graversen11}.  

This article is organized as follows: in Section~\ref{prelim}, we first recall the notions of Lévy bases, Lévy sheets and their integration theory to the extent necessary for this paper. At the beginning of Section~\ref{SPDEsec} we derive a system of SPDEs (cf. \eqref{SPDE}), which lays the groundwork for extending the CARMA process. Afterwards, we define the causal CARMA$(p,q)$ random field and the more general causal GCARMA random field, for which we drop the assumption that $A_1,...,A_d$ are in companion form. In Theorem~\ref{mildsol} we show that these random fields exist under mild assumptions and solve the SPDE system \eqref{SPDE} in the mild sense. Furthermore, we investigate the multi-parameter \emph{CARMA kernel} in more detail and present several alternative representations. This section concludes with a remark on the connection to the VOU process studied in \cite{Pham18}. Section~\ref{properties} is devoted to distributional and path properties of causal CARMA random fields. Expressions for the autocovariance function (cf. Theorem~\ref{2ndorder} and Proposition~\ref{2ndorder2}) and the spectral density (cf. Corollary~\ref{spectralden}) are derived and Theorem~\ref{pathprop} establishes some path properties.
Finally, we investigate sampling properties of causal CARMA random fields. Under a mild spectral condition, Theorem~\ref{CARMAARMA} shows that sampling on an equidistant lattice leads to a spatial ARMA process, which generally has infinitely many moving average terms. By contrast, Example~\ref{(1,1)-dependent 2} depicts a case with finitely many moving average terms.

The following notation will be used throughout this article: $C$ denotes a generic strictly positive constant which may change its value from line to line without affecting any argumentation. We use $\bone_{\{\cdot\}}$ for the indicator function so that the Heaviside function may be written as $\bone_{\{t\geq 0\}}$. If $A$ is a matrix (or a vector), then $A^\top$ denotes the transpose of $A$. The prime symbol $'$ stands for differentiation of a univariate function. For multivariate functions, we use $\partial_z$ for partial differentiation with respect to the variable $z$, or $\partial_1$ for partial differentiation with respect to the first variable. Components of a $d$-dimensional vector $u$ are denoted by $u_1,...,u_d$ if not stated otherwise. Furthermore, $\|u\|$ is the Euclidean norm, $u\cdot v\in\bbr$ is the scalar product, $u\odot v\in\bbr^d$ is the componentwise product and we write $u\leq v$ if and only if $u_i\leq v_i$ for all for $u,v\in\bbr^d$ and $i\in\{1,...,d\}$. The imaginary unit is $\bi$ and we set $\bbr_+=[0,\infty)$.

\section{Lévy bases and Lévy sheets}\label{prelim}

Throughout this article we will use homogeneous Lévy bases and we define them directly through their Lévy-Itô decomposition as a sum of a deterministic drift part, a Gaussian part, a compensated small jumps part and a large jumps part. From now on, all stochastic objects live on a fixed complete probability space $(\Omega,\calf,\bbp)$.
\begin{Definition}~
\begin{enumerate}
\item A \emph{homogeneous Lévy basis} $\La$ on $\bbr^d$ is a family of random variables indexed by the bounded Borel subsets of $\bbr^d$ such that for all $A\in\calb_\bbb(\bbr^d)$ we have
\begin{align*}
\La(A)=&~\beta \Leb_{\bbr^d}(A) + \sigma W(A) + \int_{\bbr^d}\int_\bbr \bone_A(s)z\bone_{\{|z|\leq 1\}} \,(\frakp-\frakq)(\dd s,\dd z)\nonumber \\
&~+ \int_{\bbr^d}\int_\bbr \bone_A(s)z\bone_{\{|z|> 1\}} \,\frakp(\dd s,\dd z),
\end{align*}
where
\begin{itemize}
\item $\Leb_{\bbr^d}$ is the Lebesgue measure on $\bbr^d$ and $\beta\in\bbr$, $\si\in\bbr_+$ are constants,
\item $W$ is Gaussian white noise on $\bbr^d$ such that $\var(W(A))=\Leb_{\bbr^d}(A)$ (for more details see e.g. Chapter~I in \citet{Walsh86}),
\item $\frakp$ is a Poisson random measure on $\bbr^d\times \bbr$ with intensity measure $\frakq = \Leb_{\bbr^d}\otimes\nu$, where $\nu$ is a Lévy measure on $\bbr$ (see e.g. Chapter~II in \citet{Jacod03} for more details on Poisson random measures and their integration theory).
\end{itemize}
\item The triplet $(\beta,\si^2,\nu)$ is called the \emph{characteristics} of $\La$. If $\int_\bbr |z|\bone_{\{|z| >1 \}}\,\nu(\dd z)<\infty$, we say that $\La$ has a finite first moment and define $\kappa_1:=\beta+\int_\bbr z\bone_{\{|z| >1 \}}\,\nu(\dd z)$ as the \emph{mean} of $\La$. Likewise, we say that $\La$ has a finite second moment and define $\kappa_2:=\si^2 + \int_\bbr z^2\,\nu(\dd z)$ as the \emph{variance} of $\La$ if $\int_\bbr z^2\,\nu(\dd z)<\infty$.
The \emph{cumulant generating function} (or \emph{Lévy symbol}) $\zeta\colon\bbr\to\bbc$ of $\La$ is given by
\begin{equation*}
\zeta(u)=\bi u b-\frac{1}{2}u^2\sigma^2+\int_\bbr(e^{\bi uz}-1-\bi u\bone_{\{ |z|\leq 1 \} })\nu(\dd z),\quad u\in\bbr.
\end{equation*}
\item We associate with each homogeneous Lévy basis $\La$ on $\bbr^d$ a \emph{Lévy sheet} $(L(t))_{t\in\bbr^d}$ via the equation
\begin{equation*}
L(t):=\La(\{(s_1t_1,...,s_dt_d)^\top\colon s_1,...,s_d\in[0,1]\}),\quad t\in\bbr^d.
\end{equation*}
\end{enumerate}
\halmos
\end{Definition}
The stochastic integral with respect to Lévy bases is for deterministic integrands classically defined as in \citet[Section~II]{Rajput89} (this paper uses the term infinitely divisible independently scattered random measure for Lévy basis), i.e., it is defined as the limit in probability of stochastic integrals of an approximating sequence of simple functions, where the stochastic integral for simple functions is defined canonically. We recall an integrability characterization from \cite[Theorem~2.7]{Rajput89} in the next proposition.
\begin{Proposition}\label{intcond0}
Let $g\colon\bbr^d \to\bbr$ be a measurable function and $\La$ be a homogeneous Lévy basis on $\bbr^d$ with characteristics $(\beta,\si^2,\nu)$ and cumulant generating function $\zeta$. Then the stochastic integral $\int_{\bbr^d} g(s) \,\La(\dd s)$ is well defined if and only if
\begin{enumerate}
\item $\displaystyle\int_{\bbr^d}\left| \beta g(s)+\int_\bbr (zg(s)\bone_{\{ |zg(s)|\leq 1 \}}-g(s)z\bone_{\{ |z|\leq 1 \} })\,\nu(\dd z)\right|\, \dd s<\infty$,
\item $\displaystyle\int_{\bbr^d} \sigma^2|g(s)|^2\, \dd s<\infty$,
\item $\displaystyle\int_{\bbr^d }\int_\bbr (1\wedge |zg(s)|^2)\,\nu(\dd z)\,\dd s<\infty$.
\end{enumerate}
In this case, the stochastic integral is infinitely divisible with characteristic function
\begin{align*}
\Phi\left(\int_{\bbr^d} g(s) \,\La(\dd s)\right)(u)&=\exp\left\lbrace \int_{\bbr^d} \zeta(ug(s))\,\dd s \right\rbrace\\
&=\exp\left\lbrace iu \beta_g-\frac{1}{2}u^2\sigma_g^2+\int_\bbr(e^{iuz}-1-iu\bone_{\{ |z|\leq 1 \} })\nu_g(\dd z)\right\rbrace ,\quad u\in\bbr,
\end{align*}
and characteristic triplet $(\beta_g,\sigma_g^2, \nu_g)$ given by
\begin{itemize}
\item $\beta_g=\displaystyle\int_{\bbr^d}( \beta g(s)+\int_\bbr (zg(s)\bone_{\{ |zg(s)|\leq 1 \}}-g(s)z\bone_{\{ |z|\leq 1 \} })\,\nu(\dd z))\, \dd s$,
\item $\sigma^2_g= \displaystyle\int_{\bbr^d}\sigma^2 |g(s)|^2\,\dd s$,
\item $\nu_g(B)=\displaystyle \int_{\bbr^d} \int_\bbr \bone_{\{ g(s)z\in B \}} \,\nu(\dd z)\,\dd s$ for any Borel set $B\in\calb(\bbr)$.
\end{itemize}
\end{Proposition}
The following proposition provides a sufficient integrability criterion which is easier to check and will be useful later on. For its proof, we refer to \citet{Berger17}.
\begin{Proposition}
\label{intcond} Let $g\colon\bbr^d \to\bbr$ be a measurable function such that $|g(x)|\leq C\ee^{-\eta\|x\|}$ for some positive constants $C,\eta$ and $\La$ be a homogeneous Lévy basis on $\bbr^d$ with $\int_\bbr \log(|z|)^d\bone_{\{|z|>1\}} \,\nu(\dd z)<\infty$. Then $g$ is integrable with respect to $\La$.
\end{Proposition}
Finally, the stochastic integral with respect to $L$ is defined for exactly those functions $g$ which are integrable with respect to $\La$ and we set
\begin{equation*}
\int_{\bbr^d} g(s) \,\dd L(s):=\int_{\bbr^d} g(s) \,\La(\dd s).
\end{equation*}

\section{Causal CARMA random fields as mild solutions to a system of SPDEs}\label{SPDEsec}

Our approach to defining a CARMA random field relies on a generalization of the state-space Equations \eqref{statespace1}. As the first step, we reformulate \eqref{statespace1} as
\begin{equation}\label{SDE}
\begin{aligned}
Y(t)&=b^\top X(t),\quad t\in\bbr,\\
\left(I_p\partial_t-A\right)X(t)&=c\dot L(t),\quad t\in\bbr,
\end{aligned}
\end{equation}
where $I_p$ is the identity matrix in $\bbr^{p\times p}$, $I_p\partial_t$ is a matrix whose entries are $\partial_t$ on the diagonal and zero otherwise, and $\left(I_p\partial_t-A\right)$ is a system of ordinary differential operators with constant coefficients acting on the state vector $X$. The symbol $\dot L$ denotes the formal partial differentiation of $L$ in each of $t$'s components once. In the purely temporal case, $t$ is one-dimensional and thus $\dot L(t)$ coincides with $\partial_tL(t)$. However, in $d$ dimensions we have
\begin{equation*}
\dot L(t)=\partial_1\cdots\partial_dL(t),\quad t\in\bbr^d.
\end{equation*}
In order to lift Equations \eqref{SDE} to $\bbr^d$, we iterate the system of differential operators for each of $t$'s components, that is, we consider the system of SPDEs
\begin{equation}
\label{SPDE}
\begin{aligned}
Y(t)&=b^\top X(t),\quad t\in\bbr^d,\\
\frakD_d X(t)&=c\dot L(t),\quad t\in\bbr^d,
\end{aligned}
\end{equation}
where $\frakD_d$ is the system of linear partial differential operators
\begin{equation*}
\frakD_d=(I_p\partial_d-A_d)\cdots(I_p\partial_1-A_1),
\end{equation*}
$A_i$ is the companion matrix to a monic polynomial $a_i(\cdot)$ of degree $p$ for each $i=1,...,d$ and $L$ is the Lévy sheet associated to a homogeneous Lévy basis $\La$ on $\bbr^d$. We are interested in solutions to \eqref{SPDE} and use the following notion, which is based on the random field approach of \citet{Walsh86}.
\begin{Definition}~
\begin{enumerate}
\item If $G\in (\cald'(\bbr^d))^{p\times p}$ is a matrix with entries in the space of real-valued distributions $\cald'(\bbr^d)$ and the application of the system of linear partial differential operators $\frakD_d$ on $G$ satisfies $\frakD_d G=I_p\delta_0$, where $\delta_0$ is the Dirac delta function, then $G$ is called a \emph{fundamental solution of $\frakD_d$ to the right} (see e.g. Section~3.8 in \citet{Hoermander69}).
\item If $G\in (L_{\loc}^1(\bbr^d))^{p\times p}$ is a fundamental solution of $\frakD_d$ to the right with entries in the space of locally integrable functions $L_{\loc}^1(\bbr^d)$ and the random field
\begin{equation*}
G*(c\La)(t):=\int_{\bbr^d} G(t-s)c\,\La(\dd s),\quad t\in\bbr^d,
\end{equation*}
exists in the sense of Section~\ref{prelim}, then $b^\top(G*(c\La))$ is called a \emph{mild solution to \eqref{SPDE}}.\halmos
\end{enumerate}
\end{Definition} 
\begin{Remark}
As in the purely temporal case, the derivative $\dot L$ in Equations \eqref{SPDE} does not exist in the classical sense. Nevertheless, there is a version of $L$ such that $\dot L$ exists in the distributional sense. It can then be identified with the homogeneous Lévy basis $\La$ in the sense that $\langle \dot L, \phi \rangle=\int \phi(s) \,\La(\dd s)$ for all test functions $\phi\in\cald(\bbr^d)$, where the angle brackets denote the application of the distribution $\dot L$ to $\phi$ (see Lemma~3.6 in \citet{Dalang15}). If in addition the Lévy measure $\nu$ satisfies $\int_{|x|>1} |x|^\al \,\nu(\dd x)<\infty$ for some $\al>0$, then $\dot L$ is even a random element of the space of tempered distributions $\cals'(\bbr^d)$ (see Theorem~3.13 in \cite{Dalang15}).
\halmos
\end{Remark}
Motivated by the solution formula \eqref{CARMAproc} for the CARMA process, we define the causal CARMA random field and we will subsequently show that it is a mild solution to \eqref{SPDE}.
\begin{Definition}
\label{defcarma} Let $q$ and $p$ be two non-negative integers such that $q<p$, $b=(b_0,...,b_{p-1})^\top\in\bbr^p$ with $b_q\neq 0$ and $b_i=0$ for $i>q$, $c=(0,...,0,1)^\top\in\bbr^p$, and $A_i$ be the companion matrix to a monic polynomial $a_i$ of degree $p$ with real coefficients and roots having strictly negative real parts for $i=1,...,d$. A random field $(Y(t))_{t\in\bbr^d}$ is called \emph{(causal) CARMA$(p,q)$ random field} if it satisfies the equations
\begin{equation}\label{statespace2}
\begin{aligned}
Y(t)&=b^\top X(t),\quad t\in\bbr^d,\\
X(t)&=\int_{-\infty}^{t_1}\cdots\int_{-\infty}^{t_d} \ee^{A_1(t_1-s_1)}\cdots\ee^{A_d(t_d-s_d)}c \,\La(\dd s),\quad t\in\bbr^d,
\end{aligned}
\end{equation}
where $\La$ is a homogeneous Lévy basis on $\bbr^d$ with $\int_\bbr \log(|z|)^d\bone_{\{|z|>1\}} \,\nu(\dd z)<\infty$. A (causal) CARMA$(p,0)$ random field is also called a \emph{(causal) CAR$(p)$ random field}.
\halmos
\end{Definition}
Here causality is understood in the sense that the values of $X(t)$ and $Y(t)$ at point $t\in\bbr^d$ only depend on the values of $\La$ on the set $(-\infty,t_1]\times\cdots\times(-\infty,t_d]$. Causal CARMA$(p,q)$ random fields belong to the following class, for which we drop the requirements that every $A_i$ is a companion matrix and the specific choice $c=(0,...,0,1)^\top$.
\begin{Definition}
\label{defgcarma} Let $p\geq1$ be an integer, $b,c\in\bbr^p$, $A_i\in\bbr^{p\times p}$ with eigenvalues having strictly negative real parts for $i=1,...,d$, and $\La$ be a homogeneous Lévy basis on $\bbr^d$ with $\int_\bbr \log(|z|)^d\bone_{\{|z|>1\}} \,\nu(\dd z)<\infty$. A random field $(Y(t))_{t\in\bbr^d}$ is called \emph{(causal) generalized CARMA (GCARMA) random field} if it satisfies Equations \eqref{statespace2}.
\halmos
\end{Definition}
Since each CARMA$(p,q)$ random field is also a GCARMA random field, every result which applies to GCARMA random fields also applies to CARMA$(p,q)$ random fields. On the other hand, it is easy to find GCARMA random fields which are not CARMA random fields if we fix the order $p$ (cf. for instance Example~\ref{notCARMA} in the Appendix). The next theorem shows existence of GCARMA random fields and establishes the connection to the system \eqref{SPDE}. In what follows, $\mu_i(\la_i)$ denotes the algebraic multiplicity of the eigenvalue $\la_i$ with respect to the matrix $A_i$.
\begin{Theorem}
\label{mildsol} Under the conditions of Definition~\ref{defcarma} (resp. Definition~\ref{defgcarma}) the CARMA$(p,q)$ (resp. GCARMA) random field $(Y(t))_{t\in\bbr^d}$ exists and it is a mild solution to \eqref{SPDE}.
\end{Theorem}
\begin{Proof}
For the existence we have to check that the stochastic integral in \eqref{statespace2} exists. 
Let $A_i=S_i J_i S_i^{-1}$ be a Jordan decomposition of $A_i$ for $i=1,...,d$. Then we have that $\ee^{A_i(t_i-s_i)}=S_i \ee^{J_i(t_i-s_i)} S_i^{-1}$, from which we infer that each entry of the matrix $\ee^{A_i(t_i-s_i)}$ is a (possibly complex) linear combination of $\{ (t_i-s_i)^{k_i}\ee^{\la_i(t_i-s_i)}\colon \la_i\text{ is an eigenvalue of } A_i,\, 0\leq k_i\leq \mu_i(\la_i)-1 \}$ for $i=1,...,d$. Hence, each component of the integrand $\ee^{A_1(t_1-s_1)}\cdots\ee^{A_d(t_d-s_d)}c$ is a (possibly complex) linear combination of the set
\begin{equation*}
\{ (t_1-s_1)^{k_1}\ee^{\la_1(t_1-s_1)}\cdots (t_d-s_d)^{k_d}\ee^{\la_d(t_d-s_d)}\colon \la_i\text{ is an eigenvalue of } A_i,\, 0\leq k_i\leq \mu_i(\la_i)-1 \}.
\end{equation*}
This shows that the integrability criteria of Proposition~\ref{intcond} are satisfied.

In order to show that $Y$ is a mild solution to \eqref{SPDE}, we have to show that the matrix-valued function
\begin{equation*}
G(t)=\ee^{A_1t_1}\cdots\ee^{A_dt_d}\bone_{\{t\geq0\}},\quad t\in\bbr^d,
\end{equation*}
is a fundamental solution of
\begin{equation*}
\frakD_d =(I_p\partial_d-A_d)\cdots(I_p\partial_1-A_1)
\end{equation*}
to the right. By  line~A.2.2 in the appendix of \citet{Ortner15}, an application of $\frakD_d$ on $G$ yields
\begin{align*}
&(I_p\partial_d-A_d)\cdots(I_p\partial_1-A_1)\ee^{A_1t_1}\cdots\ee^{A_dt_d}\bone_{\{t_1\geq0\}}\cdots\bone_{\{t_d\geq0\}}\\
&\quad=(I_p\partial_d-A_d)\cdots(I_p\partial_2-A_2)I_p\delta_0(t_1)\ee^{A_2t_2}\cdots\ee^{A_dt_d}\bone_{\{t_2\geq0\}}\cdots\bone_{\{t_d\geq0\}}\\
&\quad=I_p\delta_0(t_1)(I_p\partial_d-A_d)\cdots(I_p\partial_2-A_2)\ee^{A_2t_2}\cdots\ee^{A_dt_d}\bone_{\{t_2\geq0\}}\cdots\bone_{\{t_d\geq0\}}\\
&\quad=I_p\delta_0(t_1)\cdots I_p\delta_0(t_d)=I_p\delta_0(t),
\end{align*}
where the derivatives are taken in the distributional sense and we have used the tensor product of distributions in the last line (see e.g. Section~5.1 in \citet{Hoermander90}). Since $G$ is also locally integrable, this finishes the proof.
\halmos
\end{Proof}
\begin{Example}[Stable GCARMA random fields]
Let $\eta>0$, $ 0<\alpha\leq2$ and $\La$ be a symmetric $\alpha$-stable homogeneous Lévy basis with cumulant generating function
\begin{equation*}
\zeta(u)=-\eta|u|^\alpha,\quad u\in\bbr.
\end{equation*}
Since $\La$ has moments of any order strictly smaller than $\alpha$, the GCARMA random field $(Y(t))_{t\in\bbr^d}$ of Definition~\ref{defgcarma} exists and Proposition~\ref{intcond0} shows that for each $t\in\bbr^d$ the characteristic function of $Y(t)$ is
\begin{equation*}
\Phi\left( Y(t) \right)(u)=\exp\left\lbrace -\eta|u|^\alpha\int_{\bbr_+^d} |b^\top\ee^{A_1s_1}\cdots\ee^{A_ds_d}c|^\alpha\,\dd s \right\rbrace,\quad u\in\bbr.
\end{equation*}
Hence, $Y(t)$ is symmetric $\alpha$-stable with the same stability index as $\La$.
\halmos
\end{Example}
The CARMA process in Equation~\eqref{CARMAproc} is a strictly stationary process. Similarly, we have that every GCARMA random field $Y$ is strictly stationary, that is, for every $n\in\bbn$ and $\tau,t^{(1)},...,t^{(n)}\in\bbr^d$ the distributions of $(Y(t^{(1)}),...,Y(t^{(n)}))$ and $(Y(t^{(1)}+\tau),...,Y(t^{(n)}+\tau))$ are equal.
\begin{Corollary}
Suppose that $(Y(t))_{t\in\bbr^d}$ is a CARMA$(p,q)$ (resp. GCARMA) random field. Then it has the representation
\begin{equation}
\label{convo} Y(t)=(g*\La)(t):=\int_{\bbr^d} g(t-s) \,\La(\dd s),\quad t\in\bbr^d,
\end{equation}
where the \emph{CARMA$(p,q)$ kernel} (resp. \emph{GCARMA kernel}) $g$ is given by
\begin{align}
\notag g(s)&=b^\top\ee^{A_1 s_1}\cdots\ee^{A_d s_d}c\bone_{\{s\geq0\}}\\
\label{kernel} &=\sum_{\la_1}\sum_{k_1=0}^{\mu_1(\la_1)-1}\cdots\sum_{\la_d}\sum_{k_d=0}^{\mu_d(\la_d)-1} d(\la_1,k_1,...,\la_d,k_d)s_1^{k_1}\ee^{\la_1 s_1}\cdots s_d^{k_d}\ee^{\la_d s_d}\bone_{\{s\geq0\}},
\end{align}
$s=(s_1,...,s_d)\in\bbr^d$, $\{d(\la_1,k_1,...,\la_d,k_d)\}$ is a set of complex coefficients and $\sum_{\la_i}$ denotes the sum over distinct eigenvalues of $A_i$ for $i=1,...,d$. In particular, $(Y(t))_{t\in\bbr^d}$ is strictly stationary.
\end{Corollary}
\begin{Proof}
By the proof of Theorem~\ref{mildsol}, each component of $\ee^{A_1s_1}\cdots\ee^{A_ds_d}c$, and therefore also $b^\top\ee^{A_1s_1}\cdots\ee^{A_ds_d}c$, is a (possibly complex) linear combination of the set
\begin{equation*}
\{ s_1^{k_1}\ee^{\la_1s_1}\cdots s_d^{k_d}\ee^{\la_ds_d}\colon \la_i\text{ is an eigenvalue of } A_i,\, 0\leq k_i\leq \mu_i(\la_i)-1 \}.
\end{equation*}
This fact and Equations~\eqref{statespace2} imply Equations~\eqref{convo} and \eqref{kernel}.
The strict stationarity follows from Equation~\eqref{convo}.
\halmos
\end{Proof}
A direct consequence of this representation is that, under the assumption that each $A_i$ has distinct eigenvalues, $Y$ is the sum of $p^d$ dependent and possibly complex valued CAR$(1)$ random fields (cf. Proposition~2 in \cite{Brockwell11} for the temporal analog), though some of which may vanish depending on the coefficients $d(\la_1,k_1,...,\la_d,k_d)$.
Moreover, the kernel $g$ in Equation~\eqref{kernel} is \emph{anisotropic} in contrast to the isotropic CARMA random field in \eqref{isoCARMArf}. Also, $g$ is in general \emph{non-separable}, i.e., it cannot be written as a product of the form $g(s)=g_1(s_1)\cdots g_d(s_d)$.

In general, we do not have explicit formulae for the coefficients $d(\la_1,k_1,...,\la_d,k_d)$ in \eqref{kernel} since they involve the product of $d$ different matrix exponentials. However, explicit formulae can be derived in certain special cases. The next two results in Proposition~\ref{kernel4} and Theorem~\ref{kernel2} give different methods for the calculation of these coefficients provided that a CARMA$(p,q)$ random field is given.
\begin{Proposition}\label{kernel4}
Suppose that $(Y(t))_{t\in\bbr^d}$ is a CARMA$(p,q)$ random field on $\bbr^d$ such that $A_1=\cdots=A_d$ and the polynomials $a_1(\cdot)$ and $b(\cdot)$ given in Definition~\ref{defcarma} and \eqref{polynomials} have no common roots. Then its kernel $g$ as given in \eqref{kernel} can be written as
\begin{equation*}
g(s)=\sum_{\la_1} \frac{1}{(\mu_1(\la_1)-1)!}\left[ \partial_z^{\mu_1(\la_1)-1} (z-\la_1)^{\mu_1(\la_1)}\ee^{z(s_1+\cdots+s_d)}b(z)/a_1(z) \right]_{z=\la_1}\bone_{\{s\geq0\}}
\end{equation*}
for $s\in\bbr^d$. In particular, if $A_1$ has distinct eigenvalues, the CARMA$(p,q)$ kernel reduces to
\begin{equation*}
g(s)=\sum_{\la_1} \frac{b(\la_1)}{a'_1(\la_1)}\ee^{\la_1(s_1+\cdots+s_d)}\bone_{\{s\geq0\}},\quad s\in\bbr^d,
\end{equation*}
where $a'_1(\cdot)$ is the derivative of the polynomial $a_1(\cdot)$.
\end{Proposition}
\begin{Proof}
This follows directly from Lemma~2.3 in \cite{Brockwell09}.
\halmos
\end{Proof} 
\begin{Lemma}\label{residue}
Suppose that $\phi(\cdot)$ and $\theta(\cdot)$ are two (complex) polynomials such that $\phi(\cdot)$ has distinct roots, which have strictly negative real parts. Furthermore, assume that $\rho$ is a simple closed curve encircling the roots of $\phi(\cdot)$ in the complex plane. Then we have for every $s\in\bbr$ that
\begin{equation*}
\frac{1}{2\pi\bi}\int_{\rho} \frac{\theta(z)}{\phi(z)}\ee^{sz} \,\dd z=\sum_{\lambda} \frac{\theta(\lambda)}{\phi'(\lambda)}\ee^{s\lambda},
\end{equation*}
where $\sum_\lambda$ denotes the sum over the distinct roots of $\phi(\cdot)$.
\end{Lemma}
\begin{Proof}
Let $\theta(\cdot)$ has representation $\theta(z)=\sum_{k=0}^n \theta_k z^k$. Then, we observe that
\begin{equation*}
\frac{1}{2\pi\bi}\int_{\rho} \frac{\theta(z)}{\phi(z)}\ee^{sz} \,\dd z
= \sum_{k=0}^n \frac{\theta_k}{2\pi\bi}\int_{\rho} \frac{z^k}{\phi(z)}\ee^{sz} \,\dd z
= \sum_{k=0}^n \theta_k \sum_\la \frac{\la^k}{\phi'(\lambda)}\ee^{s\lambda}
= \sum_{\lambda} \frac{\theta(\lambda)}{\phi'(\lambda)}\ee^{s\lambda},
\end{equation*}
where in the second equation we have used the residue theorem and evaluated the residues.
\halmos
\end{Proof}
\begin{Theorem}\label{kernel2}~
Suppose that $(Y(t))_{t\in\bbr^d}$ is a CARMA$(p,q)$ random field such that the polynomial $a_i(\cdot)$ as given in Definition~\ref{defcarma} has distinct roots for $i=1,...,d$. Let $a_i(z)=\sum_{l=0}^{p} \al_{i,l} z^{p-l}$ and define the polynomials
\begin{equation*}
a_{i,k}(z):=\sum_{l=0}^{p-k} \al_{i,l} z^{p-k-l},
\end{equation*}
for $k=1,...,p$ and $i=1,...,d$.
Then the CARMA$(p,q)$ kernel $g$ of $Y$ as given in \eqref{kernel} can be written as
\begin{equation}\label{kernel3}
g(s)=\sum_{\la_1}\cdots\sum_{\la_d} \left( \sum_{k_1=1}^p \cdots \sum_{k_{d}=1}^p
b_{k_1-1}\frac{\la_d^{k_{d}-1}}{a_d'(\la_d)} \prod_{i=1}^{d-1} \frac{\la_i^{k_{i}-1}a_{i,k_{i+1}}(\la_i)}{a_i'(\la_i)} \right)
\ee^{\la_1 s_1+\cdots+\la_d s_d}\bone_{\{s\geq0\}},\quad s\in\bbr^d.
\end{equation}
\end{Theorem}
\begin{Proof}
Denoting the $(k,l)$-entry of the matrix $e^{A_is_i}$ with $m_{k,l}^{(i)}$, we have by the definition of the matrix product that
\begin{equation*}
b^\top\ee^{A_1 s_1}\cdots\ee^{A_d s_d}c=\sum_{k_1=1}^p \cdots \sum_{k_{d+1}=1}^p b_{k_1-1}m_{k_1,k_2}^{(1)}m_{k_2,k_3}^{(2)}\cdots m_{k_{d},k_{d+1}}^{(d)}c_{k_{d+1}},
\end{equation*}
where we use the convention that $b=(b_0,...,b_{p-1})^\top$. Theorem~2.1 in \cite{Eller87} implies that
\begin{equation*}
m_{k,l}^{(i)}=\frac{1}{2\pi\bi}\int_{\rho_i} \frac{z^{k-1}a_{i,l}(z)}{a_i(z)}\ee^{zs_i} \,\dd z,
\end{equation*}
where the contour integral is taken over a simple closed curve $\rho_i$ encircling the eigenvalues of $A_i$ in the open left half of the complex plane. As a consequence, we get that
\begin{equation}\label{help2}
g(s)=\sum_{k_1=1}^p \cdots \sum_{k_{d+1}=1}^p
b_{k_1-1} \frac{1}{2\pi\bi}\int_{\rho_1} \frac{z^{k_1-1}a_{1,k_2}(z)}{a_1(z)}\ee^{zs_1} \,\dd z \cdots \frac{1}{2\pi\bi}\int_{\rho_d} \frac{z^{k_{d}-1}a_{d,k_{d+1}}(z)}{a_d(z)}\ee^{zs_d} \,\dd z c_{k_{d+1}} \bone_{\{s\geq0\}}
\end{equation}
for $s\in\bbr^d$. Applying Lemma~\ref{residue}, we obtain that
\begin{equation*}
g(s)=\sum_{k_1=1}^p \cdots \sum_{k_{d+1}=1}^p \sum_{\la_1}\cdots\sum_{\la_d}
b_{k_1-1} \frac{\la_1^{k_1-1}a_{1,k_2}(\la_1)}{a_1'(\la_1)}\ee^{\la_1s_1}
\cdots \frac{\la_d^{k_{d}-1}a_{d,k_{d+1}}(\la_d)}{a_d'(\la_d)}\ee^{\la_ds_d} c_{k_{d+1}} \bone_{\{s\geq0\}},
\end{equation*}
which after rearranging terms and recalling that $c=(0,...,0,1)^\top$ yields Equation~\eqref{kernel3}.
\halmos
\end{Proof}
\begin{Remark}
In the setting of Theorem~\ref{kernel2}, Equation~\eqref{kernel3} reduces for $d=1$ to
\begin{equation*}
g(s)=\sum_{\la_1} \frac{b(\la_1)}{a'_1(\la_1)}\ee^{\la_1s_1}\bone_{\{s\geq0\}},\quad s\in\bbr,
\end{equation*}
which is the known kernel representation of a causal CARMA process (cf. Remark~5 in \cite{Brockwell14}). For $d=2$, Equation~\eqref{kernel3} reduces to
\begin{equation*}
g(s)=\sum_{\la_1}\sum_{\la_2} \left( \sum_{k=1}^p \frac{b(\la_1)a_{1,k}(\la_1)\la_2^{k-1}}{a'_1(\la_1)a'_2(\la_2)} \right) \ee^{\la_1 s_1+\la_2 s_2}\bone_{\{s\geq0\}},\quad s\in\bbr^2.
\end{equation*}
\halmos
\end{Remark}
\begin{Remark}
At the end of this section we reveal a connection between CARMA$(p,q)$ random fields and Volterra-type Ornstein-Uhlenbeck (VOU) processes as studied in \citet{Pham18}. A VOU process $(W(t,x))_{(t,x)\in\bbr_+\times\bbr^d}$ is a solution to the stochastic temporo-spatial integral equation
\begin{equation*}
\label{VOU} W(t,x)=V(t,x) + \int_0^t\int_{\bbr^d} W(t-s,x-y)\,\mu_{VOU}(\dd s, \dd y) + \int_0^t\int_{\bbr^d} g_{VOU}(t-s,x-y)\,\La(\dd s,\dd y),
\end{equation*}
where $\mu_{VOU}$ is a signed measure on $\bbr_+\times\bbr^d$, $g_{VOU}\colon\bbr_+\times\bbr^d\to\bbr$ is a measurable function, $V$ is a stochastic process on $\bbr_+\times\bbr^d$ and $\La$ is a homogeneous Lévy basis. Under the conditions of Theorem~3.3 in \cite{Pham18} and the particular choice of $V$ specified therein, the unique solution to this equation is given by
\begin{equation}
\label{VOUsol} W(t,x)=\int_{-\infty}^t\int_{\bbr^d} (g_{VOU}-\rho_{VOU}*g_{VOU})(t-s,x-y)\,\La(\dd s,\dd y),\quad(t,x)\in\bbr\times\bbr^d,
\end{equation}
where $*$ denotes convolution and the resolvent $\rho_{VOU}$ is another signed measure on $\bbr_+\times\bbr^d$ which is uniquely determined by $\mu_{VOU}$ through $\rho_{VOU}*\rho_{VOU}=\rho_{VOU}+\rho_{VOU}$ (cf. Proposition~2.2 in \cite{Pham18}).

We want to show that GCARMA random fields, and thus also CARMA$(p,q)$ random fields, are parametric examples of VOU processes. In order to do so, we consider a GCARMA random field $Y$ on $\bbr^{d+1}$. Since VOU processes are formulated in space and time, we write $Y$ as a function of $(t,x)=(t,x_1,...,x_d)\in\bbr^{d+1}$ instead of $t=(t_1,...,t_{d+1})\in\bbr^{d+1}$ in this remark. Further, we assume for simplicity that the matrices $A_1$,...,$A_{d+1}$ in the definition of $Y$ all have distinct eigenvalues. In this case, the GCARMA random field $Y$ satisfies
\begin{equation}
\label{VOUCARMA} Y(t,x)=\int_{-\infty}^t\cdots\int_{-\infty}^{x_d} \sum_{\la_1}\cdots\sum_{\la_{d+1}} d(\la_1,...,\la_{d+1})\ee^{\la_1(t-s)}\cdots\ee^{\la_{d+1}(x_d-y_d)} \,\La(\dd s,\dd y)
\end{equation}
for $(t,x)\in\bbr\times\bbr^d$, which can be seen from Equations \eqref{convo} and \eqref{kernel}. The task is now to find a suitable function $g_{VOU}$ and a suitable measure $\mu_{VOU}$ such that the random field $Y$ is of the form \eqref{VOUsol}. For $\mu_{VOU}$ we may choose
\begin{equation*}
\mu_{VOU}=-\la \Leb_{\bbr_+}\otimes \delta_{0,\bbr^d}
\end{equation*}
with some arbitrary real number $\la>0$. Here, $\Leb_{\bbr_+}$ denotes the Lebesgue measure on $\bbr_+$ and $\delta_{0,\bbr^d}$ is the Dirac measure on $\bbr^d$. According to Example~B.2 in \cite{Pham18}, the corresponding resolvent $\rho_{VOU}$ is then given by
\begin{equation*}
\rho_{VOU}(\dd s,\dd y) = \la\ee^{-\la s}\,\dd s \,\delta_{0,\bbr^d}(\dd y).
\end{equation*}
With this in mind, we set
\begin{equation*}
g_{VOU}(s,y)=\sum_{\la_1}\left( \frac{\la_1+\la}{\la_1}\ee^{\la_1 s}-\frac{\la}{\la_1} \right) \left( \sum_{\la_2}\cdots\sum_{\la_{d+1}} d(\la_1,...,\la_{d+1})\ee^{\la_2 y_1}\cdots\ee^{\la_{d+1}y_d}\bone_{\{y\geq 0\}} \right)
\end{equation*}
for $(s,y)\in\bbr\times\bbr^d$, and a basic calculation yields
\begin{align*}
&(\rho_{VOU}*g_{VOU})(s,y)\\
&\quad=\sum_{\la_1}\left( \frac{\la}{\la_1}\ee^{\la_1 s}-\frac{\la}{\la_1} \right) \left( \sum_{\la_2}\cdots\sum_{\la_{d+1}} d(\la_1,...,\la_{d+1})\ee^{\la_2 y_1}\cdots\ee^{\la_{d+1}y_d}\bone_{\{y\geq 0\}} \right).
\end{align*}
Plugging these two equations into \eqref{VOUsol} we observe that the two random fields in \eqref{VOUsol} and \eqref{VOUCARMA} coincide, giving us the desired result.
\halmos 
\end{Remark}

\section{Distributional and path properties}\label{properties}

In this section we examine several features of CARMA$(p,q)$ random fields. We investigate their autocovariance and their spectral density, followed by some path properties. Moreover, we analyze in detail the restriction on an equidistant discrete lattice.

\subsection{Second-order structure}

The first result in this section determines the autocovariance function. We use the convention that $\cov[V,W]$ denotes the matrix $(\cov[V_i,W_j])_{1\leq i,j\leq d}$ for any two random vectors $V,W\in\bbr^d$.
\begin{Theorem}
\label{2ndorder} Suppose that $(Y(t))_{t\in\bbr^d}$ is a CARMA$(p,q)$ (resp. GCARMA) random field.
\begin{enumerate}
\item If $\La$ has a finite first moment, then $Y(t)$ and $X(t)$ have as well for all $t\in\bbr^d$. They are given for $t\in\bbr^d$ by
\begin{equation*}
\bbe[X(t)]=\kappa_1\int_{\bbr_+^d} \ee^{A_1 s_1}\cdots\ee^{A_d s_d}c \,\dd s \quad\text{and}\quad \bbe[Y(t)]=\kappa_1\int_{\bbr_+^d} b^\top\ee^{A_1 s_1}\cdots\ee^{A_d s_d}c \,\dd s.
\end{equation*}
\item If $\La$ has further a finite second moment, then $Y(t)$ and $X(t)$ have as well for all $t\in\bbr^d$. They are given for $t\in\bbr^d$ by
\begin{equation*}
\var[X(t)]=\Sigma:=\kappa_2\int_{\bbr_+^d} \ee^{A_1 s_1}\cdots\ee^{A_d s_d}cc^\top\ee^{A^\top_d s_d}\cdots\ee^{A^\top_1 s_1} \,\dd s \quad\text{and}\quad \var[Y(t)]=b^\top\Sigma b.
\end{equation*}
In this case, the autocovariance function $\ga$ of $Y$ has the form
\begin{align}
\notag\ga(t)&=\kappa_2\sum_{\la_1}\sum_{k_1=0}^{\mu_1(\la_1)-1}\cdots\sum_{\la_d}\sum_{k_d=0}^{\mu_d(\la_d)-1}\sum_{v\in\{-1,1\}^d} d_v(\la_1,k_1,...,\la_d,k_d)\bone_{\{ t\odot v\in\bbr_+^d \}}\\
\label{auto0}
&\quad \times t_1^{k_1}\ee^{\la_1|t_1|}\cdots t_d^{k_d}\ee^{\la_d|t_d|},\quad t\in\bbr^d,
\end{align}
where $\{d_v(\la_1,k_1,...,\la_d,k_d)\}$ is a set of complex coefficients for every $v\in\{-1,1\}^d$ such that
\begin{equation*}
d_v(\la_1,k_1,...,\la_d,k_d)=d_{-v}(\la_1,k_1,...,\la_d,k_d),
\end{equation*}
$\sum_{\la_i}$ denotes the sum over distinct eigenvalues of $A_i$ for $i=1,...,d$ and $\mu_i(\la_i)$ is the algebraic multiplicity of the eigenvalue $\la_i$ with respect to the matrix $A_i$.
\item If even further all $A_i$, $i=1,...,d$, commute, then the autocovariance function $\ga$ of $Y$ has representation
\begin{equation}
\label{auto1} \ga(t)=
b^\top\ee^{A_1 |t_1|\bone_{\{t_1\geq0\}}}\cdots\ee^{A_d |t_d|\bone_{\{t_d\geq0\}}}\Sigma
\ee^{A_1^\top |t_1|\bone_{\{t_1<0\}}}\cdots\ee^{A_d^\top |t_d|\bone_{\{t_d<0\}}}b,\quad t\in\bbr^d.
\end{equation}
\end{enumerate}
\end{Theorem}
\begin{Proof}
We prove the statements for $d=2$, the proof for higher dimensions is completely analogous.
The expressions for $\bbe[X(t)]$, $\bbe[Y(t)]$, $\var[X(t)]$ and $\var[Y(t)]$ are consequences of Corollary~4.2 in \cite{Pham18}. Additionally, a consideration of the involved limits of integration yields for $r\in\bbr^2$ and $t\in\bbr^2_+$ that
\begin{equation}\label{help4}
\cov[X(r+t),X(r)]=\kappa_2\int_{\bbr_+^2} \ee^{A_1 (s_1+t_1)}\ee^{A_2 (s_2+t_2)}cc^\top\ee^{A^\top_2 s_2}\ee^{A^\top_1 s_1} \,\dd s,
\end{equation}
and
\begin{equation}\label{help5}
\cov[X(r+(t_1,-t_2)^\top),X(r)]=\kappa_2\int_{\bbr_+^2} \ee^{A_1 (s_1+t_1)}\ee^{A_2 s_2}cc^\top\ee^{A^\top_2 (s_2+t_2)}\ee^{A^\top_1 s_1} \,\dd s.
\end{equation}
Since $Y(r)=b^\top X(r)$, we have that $\cov[Y(r+t),Y(r)]=b^\top\cov[X(r+t),X(r)]b$ and $\cov[Y(r+(t_1,-t_2)^\top),Y(r)]=b^\top\cov[X(r+(t_1,-t_2)^\top),X(r)]b$. Using the symmetry of $\ga$ and a similar argument as in the proof of Theorem~\ref{mildsol}, we obtain formula \eqref{auto0}.
If the matrices $A_1$ and $A_2$ commute, then Equations \eqref{help4} and \eqref{help5} simplify to
\begin{equation*}
\cov[X(r+t),X(r)]=\ee^{A_1 t_1}\ee^{A_2 t_2}\Sigma,
\end{equation*}
and
\begin{equation*}
\cov[X(r+(t_1,-t_2)^\top),X(r)]=\ee^{A_1 t_1}\Sigma\ee^{A^\top_2 t_2},
\end{equation*}
which proves Equation~\eqref{auto1}.
\halmos
\end{Proof}
\begin{Remark}
\begin{enumerate}
\item For $d=2$ Equations~\eqref{auto0} and \eqref{auto1} simplify to
\begin{align}
\notag\ga(t)&=\kappa_2\sum_{\la_1}\sum_{k_1=0}^{\mu_1(\la_1)-1}\sum_{\la_2}\sum_{k_2=0}^{\mu_2(\la_2)-1} \Big( d_{(1,1)}(\la_1,k_1,\la_2,k_2)\bone_{\{ t_1t_2\geq0 \}}\\
\notag&\quad+d_{(1,-1)}(\la_1,k_1,\la_2,k_2)\bone_{\{ t_1t_2<0 \}} \Big) t_1^{k_1}\ee^{\la_1|t_1|}t_2^{k_2}\ee^{\la_2|t_2|},\quad t\in\bbr^2,
\end{align}
and
\begin{equation}
\notag \ga(t)=\begin{cases}
b^\top\ee^{A_1 |t_1|}\ee^{A_2 |t_2|}\Sigma b, \quad\text{if }t_1t_2\geq0,\\
b^\top\ee^{A_1 |t_1|}\Sigma\ee^{A^\top_2 |t_2|}b, \quad\text{if }t_1t_2<0.
\end{cases}
\end{equation}
\item If all matrices $A_1,...,A_d$ have distinct eigenvalues, then all eigenvalues have algebraic multiplicity one. In this case, we write $d_v(\la_1,...,\la_d)$ instead of $d_v(\la_1,0,...,\la_d,0)$ in Equation~\eqref{auto0}.
\end{enumerate}
\halmos
\end{Remark}
Theorem~\ref{2ndorder} tells us that the value of the autocovariance function $\ga(t)$ depends on the \emph{quadrant} of $\bbr^2$, and more generally on the \emph{orthant} of $\bbr^d$, in which $t$ lies. This unique second-order structure is basically induced by the causality feature of the random field $Y$. It is in particular neither isotropic, in contrast to the CARMA model in \cite{Brockwell17} (see Theorem~2 in this reference), nor separable. Also, we remark that two companion matrices commute if and only if they are equal and thus the third part of Theorem~\ref{2ndorder} actually requires $A_1=\cdots=A_d$ in the case of a CARMA$(p,q)$ random field $Y$.

The next result relates the coefficients $d_v(\la_1,k_1,...,\la_d,k_d)$ in Equation~\eqref{auto0} to the the coefficients $d(\la_1,k_1,...,\la_d,k_d)$ of the kernel $g$ in Equation~\eqref{kernel}. For brevity, we only deal with the case where each matrix $A_i$ has distinct eigenvalues and $d=2$.
\begin{Proposition}\label{2ndorder2}
Suppose that $(Y(t))_{t\in\bbr^2}$ is a CARMA$(p,q)$ (resp. GCARMA) random field on the plane $\bbr^2$ such that $\La$ has a finite second moment, both $A_1$ and $A_2$ have distinct eigenvalues and the GCARMA kernel is given by
\begin{equation*}
g(s)=\sum_{i_1,i_2=1}^p d(i_1,i_2)\ee^{\la_1(i_1)s_1}\ee^{\la_2(i_2)s_2}\bone_{\{ s\geq0 \}},\quad s\in\bbr^2,
\end{equation*}
where $(\la_n(i_n))_{1\leq i_n\leq p}$ is an enumeration of the distinct eigenvalues of $A_n$ for $n=1,2$. Then the autocovariance function of $Y$ is
\begin{equation*}
\ga(t)=\kappa_2\sum_{i_1,i_2=1}^p \left( \sum_{j_1,j_2=1}^p \frac{d(i_1,i_2)d(j_1,j_2)}{(\la_1(j_1)+\la_1(i_1))(\la_2(j_2)+\la_2(i_2))} \right) \ee^{\la_1(i_1)t_1}\ee^{\la_2(i_2)t_2},
\end{equation*}
if $t_1t_2\geq0$, and
\begin{equation*}
\ga(t)=\kappa_2\sum_{i_1,i_2=1}^p \left( \sum_{j_1,j_2=1}^p \frac{d(i_1,j_2)d(j_1,i_2)}{(\la_1(j_1)+\la_1(i_1))(\la_2(j_2)+\la_2(i_2))} \right) \ee^{\la_1(i_1)t_1}\ee^{\la_2(i_2)t_2},
\end{equation*}
if $t_1t_2<0$. In particular, we have for all $i_1,i_2=1,...,p$ that
\begin{equation*}
d_{(1,1)}(\la_1(i_1),\la_2(i_2)) = d_{(-1,-1)}(\la_1(i_1),\la_2(i_2)) = \sum_{j_1,j_2=1}^p \frac{d(i_1,i_2)d(j_1,j_2)}{(\la_1(j_1)+\la_1(i_1))(\la_2(j_2)+\la_2(i_2))},
\end{equation*}
and\begin{equation*}
d_{(1,-1)}(\la_1(i_1),\la_2(i_2)) = d_{(-1,1)}(\la_1(i_1),\la_2(i_2)) = \sum_{j_1,j_2=1}^p \frac{d(i_1,j_2)d(j_1,i_2)}{(\la_1(j_1)+\la_1(i_1))(\la_2(j_2)+\la_2(i_2))}.
\end{equation*}
\end{Proposition}
\begin{Proof}
It is sufficient to consider the case when $t_1\geq0$ and $t_2\geq0$ since all other cases follow analogously. Once again, Corollary~4.2 in \cite{Pham18} implies that
\begin{align*}
\ga(t)&=\kappa_2\int_{\bbr_+^2} g(s)g(s+t) \,\dd s\\
&=\kappa_2\int_{\bbr_+^2} \left( \sum_{i_1,i_2=1}^p d(i_1,i_2)\ee^{\la_1(i_1)s_1}\ee^{\la_2(i_2)s_2} \right)\left( \sum_{i_1,i_2=1}^p d(i_1,i_2)\ee^{\la_1(i_1)(s_1+t_1)}\ee^{\la_2(i_2)(s_2+t_2)} \right) \,\dd s\\
&=\kappa_2\sum_{i_1,i_2=1}^p d(i_1,i_2)\ee^{\la_1(i_1)t_1}\ee^{\la_2(i_2)t_2}\left( \int_{\bbr_+^2} \sum_{j_1,j_2=1}^p d(j_1,j_2)\ee^{(\la_1(j_1)+\la_1(i_1))s_1}\ee^{(\la_2(j_2)+\la_2(i_2))s_2} \,\dd s \right)\\
&=\kappa_2\sum_{i_1,i_2=1}^p d(i_1,i_2)\ee^{\la_1(i_1)t_1}\ee^{\la_2(i_2)t_2}\left( \sum_{j_1,j_2=1}^p \frac{d(j_1,j_2)}{(\la_1(j_1)+\la_1(i_1))(\la_2(j_2)+\la_2(i_2))} \right)\\
&=\kappa_2\sum_{i_1,i_2=1}^p \left( \sum_{j_1,j_2=1}^p \frac{d(i_1,i_2)d(j_1,j_2)}{(\la_1(j_1)+\la_1(i_1))(\la_2(j_2)+\la_2(i_2))} \right) \ee^{\la_1(i_1)t_1}\ee^{\la_2(i_2)t_2}.
\end{align*}
\halmos
\end{Proof}
From Theorem~\ref{2ndorder} we observe that the autocovariance function $\ga$ of a GCARMA random field is integrable over $\bbr^d$. This property is also called \emph{short-range dependency} and it implies that the spectral density $f$ of $\ga$ exists, which is defined as
\begin{equation*}
f(\om)=\frac{1}{(2\pi)^d}\int_{\bbr^d} \ga(t)\ee^{-\bi\om \cdot t} \,\dd t,\quad \om\in\bbr^d.
\end{equation*}
We present an explicit formula in the case when all matrices $A_i$ have distinct eigenvalues.
\begin{Corollary}\label{spectralden}
Suppose that $(Y(t))_{t\in\bbr^d}$ is a CARMA$(p,q)$ (resp. GCARMA) random field such that $\La$ has a finite second moment and $A_i$ has distinct eigenvalues for $i=1,...,d$. Then the spectral density $f$ of $Y$ has the form
\begin{equation}\label{specden}
f(\om)=\frac{\kappa_2}{(2\pi)^d}\sum_{\la_1}\cdots\sum_{\la_d}\sum_{v\in\{-1,1\}^d} \frac{d_v(\la_1,...,\la_d)}{(\bi v_1\om_1-\la_1)\cdots(\bi v_d\om_d-\la_d)},\quad \om\in\bbr^d.
\end{equation}
\end{Corollary}
\begin{Proof}
Equation~\eqref{auto0} shows that
\begin{equation*}
\ga(t)=\kappa_2\sum_{\la_1}\cdots\sum_{\la_d}\sum_{v\in\{-1,1\}^d} d_v(\la_1,...,\la_d)\bone_{\{ t\odot v\in\bbr_+^d \}} \ee^{\la_1|t_1|+\cdots+\la_d|t_d|},\quad t\in\bbr^d.
\end{equation*}
We also have that
\begin{equation*}
\int_\bbr \bone_{\{ t_i v_i\geq0 \}}\ee^{\la_i|t_i|}\ee^{-\bi\om_i t_i}\,\dd t_i=\frac{1}{\bi v_i\om_i-\la_i}
\end{equation*}
for $i=1,...,d$ and $\om_i\in\bbr$.
Hence, we get \eqref{specden} by applying the Fourier transform.
\halmos
\end{Proof}
\begin{Example}[Second-order structure of a CAR$(p)$ random field on the plane]~

\noindent Let $(Y(t))_{t\in\bbr^2}$ be a CAR$(p)$ random field on the plane $\bbr^2$ such that $\La$ has a finite second moment and both $A_1$ and $A_2$ have distinct eigenvalues. Recall from Equation~\eqref{help2} that the CAR$(p)$ kernel $g$ has the representation
\begin{equation*}
g(s)=\sum_{k=1}^p \frac{1}{(2\pi\bi)^2} \int_{\rho_1} \frac{a_{1,k}(z)}{a_1(z)}\ee^{zs_1} \,\dd z \int_{\rho_2} \frac{z^{k-1}}{a_2(z)}\ee^{zs_2} \,\dd z \bone_{\{s\geq0\}},\quad s=(s_1,s_2)\in\bbr^2,
\end{equation*}
where we use the notation of Theorem~\ref{kernel2}. Since $|a_{1,k}(z)/a_1(z)|=\calo(|z|^{-1})$ as $z\to\infty$, an application of Theorem~2.2 of Chapter VI in \cite{Lang99} yields
\begin{align*}
g(s)&=\sum_{k=1}^p \frac{1}{(2\pi)^2} \int_{\bbr} \frac{a_{1,k}(\bi\om_1)}{a_1(\bi\om_1)}\ee^{\bi\om_1s_1} \,\dd \om_1 \int_{\bbr} \frac{(\bi\om_2)^{k-1}}{a_2(\bi\om_2)}\ee^{\bi\om_2s_2} \,\dd \om_2\\
&=\frac{1}{(2\pi)^2}\int_{\bbr}\int_{\bbr} \sum_{k=1}^p \frac{a_{1,k}(\bi\om_1)(\bi\om_2)^{k-1}}{a_1(\bi\om_1)a_2(\bi\om_2)}\ee^{\bi(\om_1s_1+\om_2s_2)} \,\dd \om_1\,\dd \om_2,\quad s=(s_1,s_2)\in\bbr^2.
\end{align*}
This allows us to recognize that the Fourier transform of $g$ is equal to
\begin{equation*}
\tilde g(\om)=\sum_{k=1}^p \frac{a_{1,k}(\bi\om_1)(\bi\om_2)^{k-1}}{a_1(\bi\om_1)a_2(\bi\om_2)},\quad \om=(\om_1,\om_2)\in\bbr^2,
\end{equation*}
which immediately implies the spectral density
\begin{equation*}
f(\om)=\frac{\kappa_2}{(2\pi)^2}|\tilde g(\om_)|^2=\frac{\kappa_2}{(2\pi)^2} \left( \sum_{k=1}^p \frac{a_{1,k}(\bi\om_1)(\bi\om_2)^{k-1}}{a_1(\bi\om_1)a_2(\bi\om_2)} \right)
\left( \sum_{l=1}^p \frac{a_{1,l}(-\bi\om_1)(-\bi\om_2)^{l-1}}{a_1(-\bi\om_1)a_2(-\bi\om_2)} \right).
\end{equation*}
Furthermore, we conclude for the autocovariance function that
\begin{align*}
\ga(t)&=\int_{\bbr}\int_{\bbr} f(\om)\ee^{\bi(\om_1t_1+\om_2t_2)} \,\dd \om_1\,\dd \om_2\\
&=\frac{\kappa_2}{(2\pi)^2}\sum_{k,l=1}^p 
\int_{\bbr} \frac{a_{1,k}(\bi\om_1)a_{1,l}(-\bi\om_1)}{a_1(\bi\om_1)a_1(-\bi\om_1)}\ee^{\bi\om_1t_1} \,\dd \om_1 
\int_{\bbr} \frac{(\bi\om_2)^{k-1}(-\bi\om_2)^{l-1}}{a_2(\bi\om_2)a_2(-\bi\om_2)}\ee^{\bi\om_2t_2} \,\dd \om_2\\
&=\kappa_2\sum_{k,l=1}^p \Bigg( \sum_{\la_1} \frac{a_{1,k}(\la_1)a_{1,l}(-\la_1)\bone_{\{t_1\geq0\}}+a_{1,k}(-\la_1)a_{1,l}(\la_1)\bone_{\{t_1<0\}}}{a_1'(\la_1)a_1(-\la_1)} \ee^{\la_1|t_1|} \Bigg)\\
&\quad\times \Bigg( \sum_{\la_2} \frac{\la_2^{k-1}(-\la_2)^{l-1}\bone_{\{t_2\geq0\}} + (-\la_2)^{k-1}\la_2^{l-1}\bone_{\{t_2<0\}}}{a_2'(\la_2)a_2(-\la_2)} \ee^{\la_2|t_2|} \Bigg) \\
&=\kappa_2\sum_{\la_1}\sum_{\la_2}\Bigg[ \sum_{k,l=1}^p \frac{a_{1,k}(\la_1)a_{1,l}(-\la_1)\la_2^{k+l-2}}{a'_1(\la_1)a_1(-\la_1)a'_2(\la_2)a_2(-\la_2)}\\
&\quad\times\left((-1)^{l-1}\bone_{\{ t_1t_2\geq0 \}}+(-1)^{k-1}\bone_{\{ t_1t_2<0 \}}\right) \Bigg]\ee^{\la_1|t_1|+\la_2|t_2|},\quad t=(t_1,t_2)\in\bbr^2,
\end{align*}
where in the third equation we have used Lemma~\ref{residue} and Theorem~2.2 of Chapter VI in \cite{Lang99}. We remark that the procedure in this example cannot be extended to CARMA$(p,q)$ random fields with $q>0$ since $|b(z)a_{1,k}(z)/a_1(z)|=\calo(|z|^{-1})$ would not be satisfied for each $k=1,...,p$.
\halmos
\end{Example}

\subsection{Path properties}

Path properties for CARMA$(p,q)$ processes can easily be deduced from Equations \eqref{statespace1}. The sample paths of each CARMA$(p,q)$ process are $(p-q-2)$-times differentiable, provided $p>q+2$, or continuous, provided $p=q+2$, or càdlàg, provided $p=q+1$ (see e.g. Equation~(31) in \cite{Brockwell14} for more details). If additionally the driving noise is a Brownian motion, then even more regularity can be obtained.

For the class of spatial CARMA$(p,q)$ random fields it is harder to establish path properties since the system of SPDEs \eqref{SPDE} does not allow for the same reasoning as \eqref{statespace1} does. However, by drawing on maximal inequalities for multi-parameter martingales, we are able to prove the existence of a version having the path property in Definition~\ref{cadlag}, which also represents a possible generalization of the classical càdlàg property. As usual, we say that a process $(\tilde Y(t))_{t\in\bbr^d}$ is a version of the process $(Y(t))_{t\in\bbr^d}$ if for every $t\in\bbr^d$ the equality $\tilde Y(t)=Y(t)$ holds almost surely.
\begin{Definition}\label{cadlag}~
\begin{enumerate}
\item We write $v\nleq w$ if and only if $v_i> w_i$ for at least one $i\in\{1,...,d\}$. Also, for $v,w\in\bbr^d$ we define the interval $[v,w]:=\{s\in\bbr^d\colon v\leq s\leq w\}$, which may be empty.
\item A function $f\colon\bbr^d\to\bbr$ is  \emph{càdlàg} if for every $t\in\bbr^d$,
\begin{equation*}
\lim_{\substack{s\to t\\ s\geq t}}f(s)=f(t) \quad\text{and}\quad \lim_{\substack{s\to t\\ s \ngeq t}}f(s) \quad\text{exists.}
\end{equation*}
\end{enumerate}
\halmos
\end{Definition}
\begin{Theorem}
\label{pathprop} Suppose that $(Y(t))_{t\in\bbr^d}$ is a CARMA$(p,q)$ (resp. GCARMA) random field.
\begin{enumerate}
\item If the homogeneous Lévy basis $\La$ is Gaussian, then $Y$ has a version which is Hölder continuous with any exponent in $(0,1/2)$.
\item If the Lévy measure $\nu$ of $\La$ satisfies $\int_{|x|>1} |x|^\al \,\nu(\dd x)<\infty$ for some $\al\in(0,1]$, then $Y$ has a càdlàg version.
\end{enumerate}
\end{Theorem}
\begin{Proof}
We only prove the assertions in two dimensions since higher dimensions can be treated completely analogously.

\vspace{0.5\baselineskip}\noindent (1) Without loss of generality, we assume that $\beta=0$, that is, $\La$ has mean zero. Recall from Equation~\eqref{auto0} that the autocovariance function is given by
\begin{align*}
\ga(t)&=\kappa_2\sum_{\la_1}\sum_{k_1=0}^{\mu_1(\la_1)-1}\sum_{\la_2}\sum_{k_2=0}^{\mu_2(\la_2)-1} \Big( d_{(1,1)}(\la_1,\mu_1,\la_2,\mu_2)\bone_{\{ t_1t_2\geq0 \}}\\
&\quad+d_{(1,-1)}(\la_1,\mu_1,\la_2,\mu_2)\bone_{\{ t_1t_2<0 \}} \Big) t_1^{k_1}\ee^{\la_1|t_1|}t_2^{k_2}\ee^{\la_2|t_2|},\quad t\in\bbr^2.
\end{align*}
This allows us to conclude that every $t\in\bbr^d$ and every $s\in\bbr^d$ with sufficiently small norm $\|s\|$ satisfy
\begin{equation*}
\bbe[|Y(t)-Y(t+s)|^2]=2(\ga(0)-\ga(s))\leq C\|s\|.
\end{equation*}
Hence, applying Kolmogorov's continuity theorem (see e.g. Theorem~3.23 of \cite{Kallenberg02}) and the fact that $Y$ is a Gaussian process finishes the proof of the first part.

\vspace{0.5\baselineskip}\noindent (2) The proof of the second part is similar to the proof of Theorem~5.5 in \cite{Pham18}. We therefore only sketch the main ideas. Due to the first part above, we may assume that $\si^2=\beta=0$. The remaining compensated small jumps part and large jumps part are considered separately.

\vspace{\baselineskip}\noindent\textbf{Case 1: $\La(\dd s)=\int_\bbr z\bone_{\{|z|\leq 1 \}} \,(\frakp-\frakq)(\dd s,\dd z)$.}

\vspace{0.5\baselineskip}\noindent Referring to Equation \eqref{convo}, we observe that the claim clearly holds true for the process $g*\La_n$ with $\La_n(\dd s):=\int_\bbr z\bone_{\{1/n\leq|z|\leq 1 \}}\bone_{\{|s|\leq n \}} \,(\frakp-\frakq)(\dd s,\dd z)$ since the latter has only finitely many jumps. Consequently, we may close this case once we are able to show that $g*\La^n$ converges uniformly on compacts in probability to $0$, where $\La^n:=\La-\La_n$. By the fundamental theorem of calculus, we have for the GCARMA kernel $g$ that
\begin{align*}
g(t_1-s_1,t_2-s_2)&=g(-s_1,-s_2) + \int_0^{t_1} \partial_1 g(r_1-s_1,-s_2) \,\dd r_1\\
&\quad +\int_0^{t_2} \partial_2 g(-s_1,r_2-s_2) \,\dd r_2 + \int_0^{t_1}\int_0^{t_2} \partial_1\partial_2 g(r_1-s_1,r_2-s_2) \,\dd r_2\dd r_1.
\end{align*}
Putting this decomposition into $g*\La^n$, we end up with four processes which can be handled one by one. For instance, the last one satisfies
\begin{align*}
&\bbe\left[ \sup_{t\in[-v,v]}\left| \int_{-\infty}^{t_1}\int_{-\infty}^{t_2} \int_0^{t_1}\int_0^{t_2} \partial_1\partial_2 g(r_1-s_1,r_2-s_2) \,\dd r_2\dd r_1 \La^n(\dd s_1,\dd s_2) \right|^2 \right]\\
&\quad \leq C\int_{-v_1}^{v_1}\int_{-v_2}^{v_2}\bbe\left[ \sup_{t\in[-v,v]} \left| \int_{-\infty}^{t_1}\int_{-\infty}^{t_2} \partial_1\partial_2 g(r_1-s_1,r_2-s_2) \,\La^n(\dd s_1,\dd s_2) \right|^2 \right]\,\dd r_2\dd r_1\\
&\quad \leq C\int_{-v_1}^{v_1}\int_{-v_2}^{v_2}\bbe\left[ \left| \int_{-\infty}^{v_1}\int_{-\infty}^{v_2} \partial_1\partial_2 g(r_1-s_1,r_2-s_2) \,\La^n(\dd s_1,\dd s_2) \right|^2 \right]\,\dd r_2\dd r_1
\end{align*}
for some $v\in\bbr_+^2$. Note that we have used a stochastic Fubini theorem (see e.g. Theorem~2 in \cite{Lebedev96}) in the second line and Cairoli's maximal inequality (see e.g. Corollary~2.3.1 of Chapter~7 in \cite{Khoshnevisan02}) in the third line, which in turn converges to zero as $n$ tends to infinity due to the dominated convergence theorem. The three other parts in the decomposition can be dealt with analogously.

\vspace{\baselineskip}\noindent\textbf{Case 2: $\La(\dd s)=\int_\bbr z\bone_{\{|z|> 1 \}} \,\frakp(\dd s,\dd z)$.}

\vspace{0.5\baselineskip}\noindent The difference between this case and the previous one is that $\La_n(\dd s):=\int_\bbr z\bone_{\{|z|>1 \}}\bone_{\{|s|\leq n \}} \,\frakp(\dd s,\dd z)$ and instead of decomposing $g$, we directly estimate
\begin{align*}
&\bbe\left[ \sup_{t\in[-v,v]}\left| \int_{-\infty}^{t_1}\int_{-\infty}^{t_2} g(t_1-s_1,t_2-s_2) \,\La^n(\dd s_1,\dd s_2) \right|^\al \right]\\
&\quad\leq\bbe\left[ \sup_{t\in[-v,v]}\int_{-\infty}^{t_1}\int_{-\infty}^{t_2}\int_\bbr g(t_1-s_1,t_2-s_2)^\al|z|^\al\bone_{\{|z|>1 \}}\bone_{\{|s|> n \}} \,\frakp(\dd s,\dd z) \right]\\
&\quad\leq \int_{-\infty}^{v_1}\int_{-\infty}^{v_2}\sup_{t\in[-v,v]} g(t_1-s_1,t_2-s_2)^\al\bone_{\{|s|> n \}} \,\dd s_1\dd s_2\int_\bbr |z|^\al\bone_{\{|z|>1 \}} \,\dd z.
\end{align*}
The first integral in the last line is well defined and converges to zero as $n\to\infty$ by dominated convergence.
\halmos
\end{Proof}
\begin{Remark}
The notion of càdlàg functions in Definition~\ref{cadlag} is slightly stronger than the notion of \emph{lamp} functions (for \emph{limits along monotone paths}), which is for instance defined in \cite{Straf72} and also in \cite{Dalang15}.
\halmos
\end{Remark}

\subsection{Sampling on an equidistant lattice}

Real-life phenomena and data thereof are usually observed and digitally stored only for a set of discrete points even if the underlying dynamics are of a continuous nature. Therefore, it is desirable to understand the behavior of a continuous model when it is discretely sampled. If the driving Lévy process has a finite second moment, it is known that an equidistantly sampled CARMA$(p,q)$ process is always an ARMA$(p,p-1)$ process driven by a weak white noise.  This fact follows from Lemma~2.1 in \cite{Brockwell09} in conjunction with Proposition 3.2.1 in \cite{Brockwell91}. We are going to generalize these two results to higher dimensions. All results in this subsection are formulated on the plane for simplicity and we use subscripts to indicate discrete parameters.
\begin{Definition}\label{DefARMA}~
\begin{enumerate}
\item A random field $(Y_t)_{t\in\bbz^2}$ is called \emph{weakly stationary} if it has finite second moments and $\cov[Y_t,Y_s]=\cov[Y_{t-s},Y_0]=:\ga(t-s)$ for every $t,s\in\bbz^2$. It is called a \emph{white noise} if $\ga(t)=0$ for every $0\neq t\in\bbz^2$. Furthermore, a weakly stationary random field $(Y_t)_{t\in\bbz^2}$ is called \emph{$(q_1,q_2)$-dependent} if its autocovariance function $\ga$ satisfies $\ga(t)=0$ whenever $|t_1|>q_1$ or $|t_2|>q_2$, and if there are points $u,v\in\bbz^2$ such that $|u_1|=q_1$, $|v_2|=q_2$, $\ga(u)\neq0$ and $\ga(v)\neq0$.
Its \emph{spectral density} $f$ is then defined by
\begin{equation*}
f(\om)=\frac{1}{(2\pi)^2}\sum_{t\in\bbz^2} \ga(t)\ee^{-\bi\om\cdot t},\quad \om\in[-\pi,\pi]^2.
\end{equation*}
\item Let $p_1,p_2,q_1$ and $q_2$ be non-negative integers and $(Z_t)_{t\in\bbz^2}$ be a white noise on $\bbz^2$. A random field $(Y_t)_{t\in\bbz^2}$ is called an \emph{ARMA$((p_1,p_2),(q_1,q_2))$ random field} if it satisfies the equation
\begin{equation*}
\sum_{k_1=0}^{p_1}\sum_{k_2=0}^{p_2} \phi_k Y_{t-k}=\sum_{k_1=0}^{q_1}\sum_{k_2=0}^{q_2} \theta_k Z_{t-k},\quad t\in\bbz^2,
\end{equation*}
where $\phi_k, \theta_k\in\bbc$ are coefficients such that $\phi_{\textbf{0}}\neq 0$, $\theta_{\textbf{0}}\neq 0$, at least one of the coefficients $\phi_{(p_1,\cdot)}$ is non-zero and similarly for $\phi_{(q_1,\cdot)}$,$\phi_{(\cdot,p_2)}$ and $\phi_{(\cdot,q_2)}$. This random field is also called \emph{AR$((p_1,p_2))$ random field} if $q_1=q_2=0$ and \emph{MA$((q_1,q_2))$ random field} if $p_1=p_2=0$.
\item Let $P$ and $Q$ be two non-empty subsets of $\bbz^2$ and $(Z_t)_{t\in\bbz^2}$ be a white noise on $\bbz^2$. A random field $(Y_t)_{t\in\bbz^2}$ is called an \emph{ARMA$(P,Q)$ random field} if it satisfies the equation
\begin{equation*}
\sum_{k\in P} \phi_k Y_{t-k}=\sum_{k\in S} \theta_k Z_{t-k},\quad t\in\bbz^2,
\end{equation*}
where $\phi_k, \theta_k\in\bbc$ are non-zero coefficients. This random field is also called \emph{AR$(P)$ random field} if $Q=\{ (0,0) \}$ and \emph{MA$(Q)$ random field} if $P=\{ (0,0) \}$.
\halmos
\end{enumerate}
\end{Definition}
Before we investigate the general GCARMA random field, let us consider the special case of a CAR$(1)$ random field $(Y(t))_{t\in\bbr^2}$ first. Assuming that $b=1$, we have the representation
\begin{equation*}
Y(t)=\int_{-\infty}^{t_1}\int_{-\infty}^{t_2} \ee^{\la_1(t_1-s_1)+\la_2(t_2-s_2)} \,\La(\dd s),\quad t\in\bbr^2,
\end{equation*}
where the real numbers $\la_1$ and $\la_2$ are strictly negative. This allows us to observe that
\begin{align}
\notag Y(t_1,t_2)&=e^{\la_1}Y(t_1-1,t_2)+e^{\la_2}Y(t_1,t_2-1)-e^{\la_1+\la_2}Y(t_1-1,t_2-1)\\
\label{help3} &\quad+\int_{t_1-1}^{t_1}\int_{t_2-1}^{t_2} \ee^{\la_1(t_1-s_1)+\la_2(t_2-s_2)} \,\La(\dd s),\quad t\in\bbr^2.
\end{align}
Setting $Y_t:=Y(t)$ for $t\in\bbz^2$, we conclude that the sampled random field $(Y_t)_{t\in\bbz^2}$ is an AR$((1,1))$ random field driven by the i.i.d. noise $Z_t:=\int_{t_1-1}^{t_1}\int_{t_2-1}^{t_2} \ee^{\la_1(t_1-s_1)+\la_2(t_2-s_2)} \,\La(\dd s)$. With a little more effort, this procedure carries over to a more general case.
\begin{Proposition}\label{sampling}
Suppose that $(Y(t))_{t\in\bbr^2}$ is a GCARMA random field on the plane $\bbr^2$ such that $A_1$ and $A_2$ commute. Then the sampled random field $(Y_t)_{t\in\bbz^2}$ satisfies the equation
\begin{equation}\label{ARMA}
\sum_{k_1,k_2=0}^{p} d_k Y_{t-k}=\sum_{k_1,k_2=0}^{p-1} \Bigg(\sum_{l_1=0}^{k_1}\sum_{l_2=0}^{k_2} d_lb^\top\ee^{(k_1-l_1)A_1+(k_2-l_2)A_2}\Bigg) R_{t-k},\quad t\in\bbz^2,
\end{equation}
where the coefficients $d_{k}\in\bbc$ are given by
\begin{align*}
&d_{0,0}z^p+d_{1,0}z^{p-1}+\cdots+d_{p,0}=\chi_{e^{A_1}}(z):=\prod_{\la_1}\prod_{k_1=0}^{\mu_1(\la_1)-1} (z-e^{\la_1})\\
&d_{0,0}z^p+d_{0,1}z^{p-1}+\cdots+d_{0,p}=\chi_{e^{A_2}}(z):=\prod_{\la_2}\prod_{k_2=0}^{\mu_2(\la_2)-1} (z-e^{\la_1})\\
&d_{k_1,k_2}=d_{k_1,0}d_{0,k_2},\quad k_1,k_2=1,...,p,
\end{align*}
and the multivariate i.i.d. noise $(R_t)_{t\in\bbz^2}$ is given by
\begin{equation*}
R_t=\int_{t_1-1}^{t_1}\int_{t_2-1}^{t_2} \ee^{A_1(t_1-s_1)+A_2(t_2-s_2)}c \,\La(\dd s),\quad t\in\bbz^2.
\end{equation*}
In particular, the right-hand side of \eqref{ARMA} is a $(p-1,p-1)$-dependent random field if $\La$ has a finite second moment.
\end{Proposition}
\begin{Proof}
We extend the proof of Lemma~2.1 in \citet{Brockwell09} which requires several additional steps that do not appear in the one-dimensional case. To this end, we show by induction that for all $t\in\bbz^2$, $n\in\bbn\cup\{0\}$ and coefficients $f_k\in\bbc$ with $k_1,k_2=0,1,...,n$ we have
\begin{align}\label{help1}
\sum_{k_1,k_2=0}^{n} f_kX_{t-k}&=\sum_{k_1,k_2=0}^{n-1} \Bigg( \sum_{l_1=0}^{k_1}\sum_{l_2=0}^{k_2} f_l\ee^{(k_1-l_1)A_1+(k_2-l_2)A_2} \Bigg) R_{t-k}\\
\notag&\quad+\sum_{k_2=0}^{n-1}\Bigg( \sum_{k_1=0}^{n} f_k\ee^{(n-k_1)A_1} \Bigg) X_{t_1-n,t_2-k_2}+\sum_{k_1=0}^{n-1}\Bigg( \sum_{k_2=0}^{n} f_k\ee^{(n-k_2)A_2} \Bigg) X_{t_1-k_1,t_2-n}\\
\notag&\quad+\Bigg( f_{n,n}I_p-\sum_{k_1,k_2=0}^{n-1} f_k\ee^{(n-k_1)A_1+(n-k_2)A_2} \Bigg) X_{t_1-n,t_2-n}
=:S_1+S_2+S_3+S_4.
\end{align}
The case $n=0$ is trivial. Assuming that the statement is valid for some $n$, we observe that
\begin{align*}
S_2&=\Bigg( \sum_{k_1=0}^{n} f_{k_1,0}\ee^{(n-k_1)A_1} \Bigg)\Bigg( e^{A_1}X_{t_1-n-1,t_2}+e^{A_2}X_{t_1-n,t_2-1}+e^{A_1+A_2}X_{t_1-n-1,t_2-1}+R_{t_1-n,t_2} \Bigg)\\
&\quad +\sum_{k_2=1}^{n-1}\Bigg( \sum_{k_1=0}^{n} f_k\ee^{(n-k_1)A_1} \Bigg) X_{t_1-n,t_2-k_2}\\
&=\Bigg( \sum_{k_1=0}^{n} f_{k_1,0}\ee^{(n+1-k_1)A_1} \Bigg)X_{t_1-n-1,t_2}
+\Bigg( \sum_{k_1=0}^{n} f_{k_1,0}\ee^{(n-k_1)A_1} \Bigg)R_{t_1-n,t_2}\\
&\quad -\Bigg( \sum_{k_1=0}^{n} f_{k_1,0}\ee^{(n+1-k_1)A_1+A_2} \Bigg)X_{t_1-n-1,t_2-1}
+\Bigg( \sum_{k_1=0}^{n} f_{k_1,0}\ee^{(n-k_1)A_1+A_2} \Bigg)X_{t_1-n,t_2-1}\\
&\quad +\sum_{k_2=1}^{n-1}\Bigg( \sum_{k_1=0}^{n} f_k\ee^{(n-k_1)A_1} \Bigg) X_{t_1-n,t_2-k_2}\\
&=\sum_{k_2=0}^{M} \Bigg( \sum_{k_1=0}^{n} f_k\ee^{(n+1-k_1)A_1} \Bigg)X_{t_1-n-1,t_2-k_2}
+\sum_{k_2=0}^{M} \Bigg( \sum_{l_1=0}^{n}\sum_{l_2=0}^{k_2} f_l\ee^{(n-l_1)A_1+(k_2-l_2)A_2} \Bigg)R_{t_1-n,t_2-k_2}\\
&\quad-\Bigg( \sum_{k_2=0}^{M}\sum_{k_1=0}^{n} f_k\ee^{(n+1-k_1)A_1+(M+1-k_2)A_2} \Bigg)X_{t_1-n-1,t_2-M-1}\\
&\quad+\Bigg( \sum_{k_2=0}^{M}\sum_{k_1=0}^{n} f_k\ee^{(n-k_1)A_1+(M+1-k_2)A_2} \Bigg)X_{t_1-n,t_2-M-1}
+\sum_{k_2=M+1}^{n-1}\Bigg( \sum_{k_1=0}^{n} f_k\ee^{(n-k_1)A_1} \Bigg) X_{t_1-n,t_2-k_2}\\
&=\sum_{k_2=0}^{n-1} \Bigg( \sum_{k_1=0}^{n} f_k\ee^{(n+1-k_1)A_1} \Bigg)X_{t_1-n-1,t_2-k_2}
+\sum_{k_2=0}^{n-1} \Bigg( \sum_{l_1=0}^{n}\sum_{l_2=0}^{k_2} f_l\ee^{(n-l_1)A_1+(k_2-l_2)A_2} \Bigg)R_{t_1-n,t_2-k_2}\\
&\quad-\Bigg( \sum_{k_2=0}^{n-1}\sum_{k_1=0}^{n} f_k\ee^{(n+1-k_1)A_1+(n-k_2)A_2} \Bigg)X_{t_1-n-1,t_2-n}\\
&\quad+\Bigg( \sum_{k_2=0}^{n-1}\sum_{k_1=0}^{n} f_k\ee^{(n-k_1)A_1+(n-k_2)A_2} \Bigg)X_{t_1-n,t_2-n}=:U_1+U_2+U_3+U_4,
\end{align*}
where $M$ is just an induction parameter ranging from $0$ to $n-1$ and we have used a similar calculation as in \eqref{help3}. By symmetry we have that
\begin{align*}
S_3&=\sum_{k_1=0}^{n-1} \Bigg( \sum_{k_2=0}^{n} f_k\ee^{(n+1-k_2)A_2} \Bigg)X_{t_1-k_1,t_2-n-1}
+\sum_{k_1=0}^{n-1} \Bigg( \sum_{l_1=0}^{k_1}\sum_{l_2=0}^{n} f_l\ee^{(k_1-l_1)A_1+(n-l_2)A_2} \Bigg)R_{t_1-k_1,t_2-n}\\
&\quad-\Bigg( \sum_{k_2=0}^{n}\sum_{k_1=0}^{n-1} f_k\ee^{(n-k_1)A_1+(n+1-k_2)A_2} \Bigg)X_{t_1-n,t_2-n-1}\\
&\quad+\Bigg( \sum_{k_2=0}^{n}\sum_{k_1=0}^{n-1} f_k\ee^{(n-k_1)A_1+(n-k_2)A_2} \Bigg)X_{t_1-n,t_2-n}=:V_1+V_2+V_3+V_4.
\end{align*}
Moreover, we may sum up
\begin{align*}
S_4+U_4+V_4&=\Bigg( \sum_{k_1,k_2=0}^{n} f_k\ee^{(n-k_1)A_1+(n-k_2)A_2} \Bigg)\Bigg( e^{A_1}X_{t_1-n-1,t_2-n}+e^{A_2}X_{t_1-n,t_2-n-1}\\
&\quad -e^{A_1+A_2}X_{t_1-n-1,t_2-n-1}+R_{t_1-n,t_2-n} \Bigg)=:W_1+W_2+W_3+W_4,
\end{align*}
and also
\begin{align*}
\sum_{k_1,k_2=0}^{n} f_kX_{t-k}&=S_1+S_2+S_3+S_4\\
&=(S_1+U_2+V_2+W_4)+(U_1+U_3+W_1)+(V_1+V_3+W_2)+W_3\\
&=\sum_{k_1,k_2=0}^{n} \Bigg( \sum_{l_1=0}^{k_1}\sum_{l_2=0}^{k_2} f_l\ee^{(k_1-l_1)A_1+(k_2-l_2)A_2} \Bigg) R_{t-k}\\
\notag&\quad+\sum_{k_2=0}^{n}\Bigg( \sum_{k_1=0}^{n} f_k\ee^{(n+1-k_1)A_1} \Bigg) X_{t_1-n-1,t_2-k_2}
+\sum_{k_1=0}^{n}\Bigg( \sum_{k_2=0}^{n} f_k\ee^{(n+1-k_2)A_2} \Bigg) X_{t_1-k_1,t_2-n-1}\\
\notag&\quad-\sum_{k_1,k_2=0}^{n} f_k\ee^{(n+1-k_1)A_1+(n+1-k_2)A_2}X_{t_1-n-1,t_2-n-1},
\end{align*}
which is equivalent to equation \eqref{help1} with $n+1$ instead of $n$. By choosing $n=p$ and $f_k=d_k$ for $k_1,k_2=0,1,...,p$, the Cayley-Hamilton theorem implies that $S_2$ and $S_3$ in \eqref{help1} vanish since $\chi_{e^{A_1}}(z)$ and $\chi_{e^{A_2}}(z)$ are the characteristic polynomials of $e^{A_1}$ and $e^{A_2}$, respectively. Finally, we have that
\begin{align*}
d_{p,p}I_p-\sum_{k_1,k_2=0}^{p-1} d_k\ee^{(p-k_1)A_1+(p-k_2)A_2}&=d_{n,n}I_p-\Bigg( \sum_{k_1=0}^{p-1} d_{k_1,0}\ee^{(p-k_1)A_1} \Bigg)\Bigg( \sum_{k_2=0}^{p-1} d_{0,k_2}\ee^{(p-k_2)A_2} \Bigg)\\
&=d_{p,p}I_p-(-d_{p,0}I_p)(-d_{0,p}I_p)=0,
\end{align*}
which shows that $S_4$ vanishes, too. By multiplying $b^\top$ to the left of \eqref{help1}, we arrive at equation \eqref{ARMA}.
\halmos
\end{Proof}
Equation~\eqref{ARMA} implies that $(Y_t)_{t\in\bbz^2}$ satisfies an autoregression of order $(p,p)$ driven by a $(p-1,p-1)$-dependent noise.
Proposition 3.2.1 in \cite{Brockwell91} states that a stationary time-discrete $q$-dependent process is a moving average process of order $q$, which in turn is established by projecting the process into the past with respect to the natural order of time in order to create the white noise sequence. However, projecting on the past with respect to the partial order $\leq$ on $\bbz^2$ does not necessarily lead to spatial white noise since this order is only a partial order and not a total order. By contrast, the lexicographic order is total and allows us to generalize Proposition 3.2.1 in \cite{Brockwell91}.
\begin{Definition}\label{lexico}
For $v,w\in\bbz^2$ we define $v\preceq w$ if and only if $v_1 = w_1$ and $v_2 \leq w_2$ or $v_1 < w_1$. Furthermore, we write $[v,w]_\preceq:=\{s\in\bbz^2\colon v\preceq s\preceq w\}$, which might be empty.
\halmos
\end{Definition}
\begin{Proposition}\label{q-dependency}
Let $(Y_t)_{t\in\bbz^2}$ be a weakly stationary $(q_1,q_2)$-dependent random field.
If its spectral density $f$ satisfies $\log f\in L^1([-\pi,\pi]^2)$, then $Y$ is a MA$([(0,0),(q_1,q_2)]_\preceq)$ random field.
\end{Proposition}
\begin{Proof}
First of all, the $(q_1,q_2)$-dependency of $Y$ implies that its spectral measure is absolutely continuous. Since $\log f\in L^1([-\pi,\pi]^2)$, Theorems~1.1.2 and 1.1.4 in \citet{Korezlioglu86} imply that $Y$ satisfies the Wold decomposition
\begin{equation*}
Y_t=\sum_{k\succeq(0,0)} \theta_k Z_{t-k},\quad t\in\bbz^2,
\end{equation*}
where $\theta_k\in\bbc$ are such that $\theta_{(0,0)}\neq0$ and $\sum_{k\succeq(0,0)}|\theta_k|^2<\infty$ and the white noise $Z$ is given by 
\begin{equation*}
Z_t=Y_t-(Y_t/H^{1+}_{t_1-1,t_2-1}),\quad t\in\bbz^2.
\end{equation*}
Here $(Y_t/H^{1+}_{t_1-1,t_2-1})$ denotes the orthogonal projection of $Y_t$ on the closed linear subspace $H^{1+}_{t_1-1,t_2-1}$ of the Hilbert space $L^2(\Omega,\calf,\bbp)$, which is generated by $\{ Y_\bs\colon s_1<t_1,s_2\in\bbz \}$ and $\{ Y_\bs\colon s_1=t_1,s_2<t_2 \}$.
Exploiting the $(q_1,q_2)$-dependency once again, we see that actually
\begin{equation*}
Y_t=\sum_{(0,0)\preceq k\preceq(q_1,q_2)} \theta_k Z_{t-k},\quad t\in\bbz^2.
\end{equation*}
\halmos
\end{Proof}
\begin{Theorem}\label{CARMAARMA}
Suppose that $(Y(t))_{t\in\bbr^2}$ is a GCARMA random field on the plane $\bbr^2$ such that $A_1$ and $A_2$ commute and the spectral density $f$ of the right-hand side of \eqref{ARMA} satisfies $\log f\in L^1([-\pi,\pi]^2)$. Then the sampled random field $(Y_t)_{t\in\bbz^2}$ is an ARMA$([(0,0),(p,p)],[(0,0),(p-1,p-1)]_\preceq)$ random field.

Furthermore, the driving spatial white noise of this ARMA random field is i.i.d. noise in each of the following cases:
\begin{itemize}
\item $p=1$.
\item $b^\top$ is a common left eigenvector of both $A_1$ and $A_2$.
\item $\La$ is Gaussian.
\end{itemize}
\end{Theorem}
\begin{Proof}
The first part follows from Proposition~\ref{sampling} and Proposition~\ref{q-dependency}.
As for the second part, we have seen in Equation~\eqref{help3} that the driving noise is i.i.d. for $p=1$. If $b^\top$ is a common left eigenvector of both $A_1$ and $A_2$, then the GCARMA random field $(Y(t))_{t\in\bbr^2}$ reduces to a CAR$(1)$ random field. In the Gaussian case we have that every white noise is actually i.i.d. noise.

\halmos
\end{Proof}
\begin{Remark}
\begin{enumerate}
\item The order defined in Definition~\ref{lexico} is more precisely called the column-by-column lexicographic order. By symmetry, Theorem~\ref{CARMAARMA} also holds for the row-by-row lexicographic order correspondingly.
\item In respect of the second part of Theorem~\ref{CARMAARMA} we note that if both $A_1$ and $A_2$ have distinct eigenvalues, then they have the same left eigenvectors since we have assumed that they commute.
\item If $A_1$ is a companion matrix, then the vector $v=(v_0,...,v_{p-1})$ is a left eigenvector of $A_1$ to the eigenvalue $\la_1$ if and only if the polynomial $v(z):=v_{p-1}z^{p-1}+\cdots+v_1z+v_0$ satisfies $v(z)=v_{p-1}a_1(z)/(z-\la_1)$, where $a_1(z)$ is the corresponding polynomial to $A_1$. In particular, $b^\top$ cannot be a left eigenvector of $A_1$ if we assume that $b(z)$ and $a_1(z)$ do not have common roots.
\end{enumerate}
\halmos
\end{Remark}
Every sampled GCARMA random field is an ARMA random field according to Theorem~\ref{CARMAARMA}. However, the MA part of this random field has infinitely many terms unless $p=1$. For instance, if we sample a CARMA$(2,1)$ random field, we obtain an ARMA$([(0,0),(2,2)],[(0,0),(1,1)]_\preceq)$ random field, where $[(0,0),(1,1)]_\preceq=\{(0,u)\in\bbz^2\colon u\geq0 \}\cup \{(1,u)\in\bbz^2\colon u\leq1 \}$.
In analogy to the purely temporal case it would be desirable to have that the $(p-1,p-1)$-dependent random field on the right-hand side of Equation~\eqref{ARMA} has a MA$(p-1,p-1)$ representation such that the sampled random field is an ARMA$((p,p),(p-1,p-1))$ random field. The next two examples illustrate that unfortunately this is not always the case.
\begin{Example}[$(1,1)$-dependent random field with no MA$(1,1)$ representation]\label{(1,1)-dependent}
Let $(Y(t))_{ t\in\bbr^2 }$ be a GCARMA random field  with parameters $b=c=(1,1)^\top$,
\begin{equation*}
A_1=A_2=\begin{pmatrix}-1 &0 \\0 &-2\end{pmatrix},
\end{equation*}
and kernel
\begin{equation*}
g(s)=b^\top\ee^{A_1 s_1}\ee^{A_2 s_2}c\bone_{\{s\geq0\}}=\left(\ee^{-(s_1+s_2)}+\ee^{-2(s_1+s_2)}\right)\bone_{\{s\geq0\}},\quad s\in\bbr^2.
\end{equation*}
Further, we assume that the variance of $\La$ satisfies $\kappa_2=1$ and denote the $(1,1)$-dependent right-hand side of Equation~\eqref{ARMA} as $(U_t)_{ t\in\bbz^2 }$. By Proposition~\ref{sampling} and straight forward calculations, the autocovariance $\hat\ga$ of $U$ satisfies
\begin{align*}
\hat\ga(0,0)&=\frac{\left(\ee^2-1\right)^2 \left(77+100 \ee^2+222 \ee^4+100 \ee^6+77 \ee^8\right)}{144 \ee^{12}},\\
\hat\ga(1,0)=\hat\ga(0,1)=\hat\ga(-1,0)=\hat\ga(0,-1)&=-\frac{\left(\ee^2-1\right)^2 \left(25+52 \ee+59 \ee^2+16 \ee^3+59 \ee^4+52 \ee^5+25 \ee^6\right)}{144 \ee^{11}},\\
\hat\ga(1,1)=\hat\ga(-1,-1)&=\frac{25+52 \ee^2-32 \ee^3-90 \ee^4-32 \ee^5+52 \ee^6+25 \ee^8}{144 \ee^{10}},\\
\hat\ga(1,-1)=\hat\ga(-1,1)&=\frac{9+32 \ee+36 \ee^2-154 \ee^4+36 \ee^6+32 \ee^7+9 \ee^8}{144 \ee^{10}}.
\end{align*}
All other values of $\hat\ga$ are zero.
Having determined the autocovariance of $U$ explicitly, we try to match $\hat\ga$ with the autocovariance of a MA$(1,1)$ random field. A generic MA$(1,1)$ random field is given by
\begin{equation*}
W_t=\sum_{k_1,k_2=0}^{1} \theta_k Z_{t-k},\quad t\in\bbz^2,
\end{equation*}
with spatial white noise $Z$ and complex coefficients $\theta_k$. Its autocovariance $\ga$ satisfies
\begin{align*}
\ga(0,0)&=|\theta_{00}|^2+|\theta_{10}|^2+|\theta_{01}|^2+|\theta_{11}|^2,\\
\ga(1,0)=\bar\ga(-1,0)&=\theta_{00}\bar\theta_{10}+\theta_{01}\bar\theta_{11},\\
\ga(0,1)=\bar\ga(0,-1)&=\theta_{10}\bar\theta_{11}+\theta_{00}\bar\theta_{01},\\
\ga(1,1)=\bar\ga(-1,-1)&=\theta_{00}\bar\theta_{11},\\
\ga(1,-1)=\bar\ga(-1,1)&=\theta_{01}\bar\theta_{10}.
\end{align*}
Again, all other values of $\ga$ are zero.
Extracting imaginary and real parts and using Gröbner bases (see e.g. Chapter~2 of \citet{Cox15} for more details) together with a computer algebra system such as Mathematica, we conclude that the system $\hat\ga=\ga$ has no complex solutions for $\{\theta_{00},\theta_{10},\theta_{01},\theta_{11} \}$. Hence, $U$ is not a MA$(1,1)$ random field.
However, the spectral density $f$ of $U$ has representation
\begin{align*}
f(\om)&=\frac{1}{(2\pi)^2}\big( \hat\ga(0,0)+2\hat\ga(1,0)\cos(\om_1)+2\hat\ga(0,1)\cos(\om_2)\\&\quad\qquad\qquad+2\hat\ga(1,1)\cos(\om_1+\om_2)+2\hat\ga(1,-1)\cos(\om_1-\om_2) \big),\quad \om\in[-\pi,\pi]^2,
\end{align*}
which is strictly positive. Consequently, $\log f$ is integrable over $[-\pi,\pi]^2$ and Theorem~\ref{CARMAARMA} yields that  $(Y_t)_{t\in\bbz^2}$ is an ARMA$([(0,0),(2,2)],[(0,0),(1,1)]_\preceq)$ random field.
\halmos
\end{Example}
\begin{Example}[$(1,1)$-dependent random field with a MA$(1,1)$ representation]\label{(1,1)-dependent 2}
We keep the setting of Example~\ref{(1,1)-dependent} with the only change that
\begin{equation*}
A_1=\begin{pmatrix}-1 &0 \\0 &-2\end{pmatrix}\quad\text{and}\quad A_2=\begin{pmatrix}-1 &0 \\0 &-1\end{pmatrix}.
\end{equation*}
In this case, the system $\hat\ga=\ga$ has eight different real solutions for $\{\theta_{00},\theta_{10},\theta_{01},\theta_{11} \}$. The exact algebraic expressions for these solutions are very lengthy and can be computed with the software Mathematica. For illustration, we present the rounded values of one of these solutions, namely
\begin{equation*}
\theta_{00}=0.752991,\quad
\theta_{10}=-0.176944,\quad
\theta_{01}=-0.277010,\quad
\theta_{11}=0.065094.
\end{equation*}
Since $U$ has the second-order structure of a MA$(1,1)$ random field, it also is a MA$(1,1)$ random field due to Theorem~10 in \citet{Karhunen47}. Therefore, $(Y_t)_{ t\in\bbz^2 }$ is indeed an ARMA$((2,2),(1,1))$ random field.
\halmos
\end{Example}
Further research has to be done to determine explicit necessary and sufficient conditions for the right-hand side of Equation~\eqref{ARMA} to be a MA$(p-1,p-1)$ random field.

\begin{appendix}
\section{Appendix}

\begin{Example}[GCARMA but not CARMA]\label{notCARMA} Let $(Y(t))_{ t\in\bbr^2 }$ be a GCARMA random field  with parameters $b=c=(1,1)^\top$,
\begin{equation*}
A_1=\begin{pmatrix}-2 &0 \\0 &-3\end{pmatrix}\quad\text{and}\quad A_2=\begin{pmatrix}-5 &0 \\0 &-7\end{pmatrix},
\end{equation*}
and kernel
\begin{equation*}
g(s)=b^\top\ee^{A_1 s_1}\ee^{A_2 s_2}c\bone_{\{s\geq0\}}=\left(\ee^{-2s_1-5s_2}+\ee^{-3s_1-7s_2}\right)\bone_{\{s\geq0\}},\quad s\in\bbr^2.
\end{equation*}
In order to check whether $Y$ has a CARMA$(2,1)$ representation, we have to find two companion matrices $\hat A_1,\hat A_2\in\bbr^{2\times2}$ and a vector $\hat b=(\hat b_0,\hat b_1)^\top\in\bbr^2$ such that
\begin{equation}\label{help7}
g(s)=\hat b^\top\ee^{\hat A_1 s_1}\ee^{\hat A_2 s_2}(0,1)^\top\bone_{\{s\geq0\}},\quad s\in\bbr^2.
\end{equation}
Observing the exponentials, we conclude that
\begin{equation*}
\hat A_1=\begin{pmatrix}0 &1 \\-6 &-5\end{pmatrix}\quad\text{and}\quad \hat A_2=\begin{pmatrix}0 &1 \\35 &-12\end{pmatrix}
\end{equation*}
have to hold. Plugging these into \eqref{help7} implies
\begin{align*}
\ee^{-2s_1-5s_2}+\ee^{-3s_1-7s_2}&=\frac{\ee^{-3s_1-7s_2}}{2}\bigg[ \hat b_0\left( -5+3\ee^{2s_2}+4\ee^{s_1}-2\ee^{s_1+2s_2} \right)\\
&\quad + \hat b_1\left( 15-9\ee^{2s_2}-8\ee^{s_1}+4\ee^{s_1+2s_2} \right) \bigg],\quad s\in\bbr^2,
\end{align*}
which has no solution for $\hat b$.
\halmos
\end{Example}

\end{appendix}
\subsection*{Acknowledgement}
The author is very much indebted to Claudia Klüppelberg for continuous support and helpful comments. Cordial thanks also go to David Berger, Carsten Chong and Alexander Lindner for inspiring discussions and helpful advice, and to Carlos Am\'endola for providing Example~\ref{(1,1)-dependent 2}. Moreover, support from the graduate program TopMath at the Technical University of Munich is acknowledged.

\addcontentsline{toc}{section}{References}

\end{document}